\documentclass[reqno,11pt]{amsart}

\usepackage{amsthm,amssymb,amsmath,amsfonts}
\usepackage{graphicx,psfrag}
\usepackage{subcaption}
\usepackage{comment}
\usepackage{enumerate}
\usepackage[utf8]{inputenc}
\usepackage[english]{babel}
\usepackage[T1]{fontenc}
\usepackage{lmodern,microtype}
\usepackage{mathpazo}
\usepackage{xcolor}
\colorlet{tn/color/defi}{red!50!black}
\usepackage{hyperref}
\hypersetup{%
    pdfauthor = {F. BOTLER, A. JIMÉNEZ, C. N. LINTZMAYER, A. PASTINE, D. A. QUIROZ, AND M. SAMBINELLI}, 
    pdftitle = {Biclique immersions in graphs with independence number 2},%
    pdfkeywords = {immersions, biclique},%
    colorlinks,%
    linkcolor = {red!60!black},%
    citecolor = {green!60!black},%
    urlcolor = {blue!60!black}%
}

\usepackage{tikz,nicefrac,pifont,fp}
\usetikzlibrary{arrows,snakes,backgrounds,patterns,matrix,shapes,fit,calc,shadows,plotmarks,arrows.meta,positioning,decorations.pathmorphing,decorations.pathreplacing,decorations.markings,fixedpointarithmetic,shapes.geometric}

\makeatletter
\let\origsection=\section 
\def\section{\@ifstar{\origsection*}{\mysection}}
\def\mysection{\@startsection{section}{1}\z@{.7\linespacing\@plus\linespacing}{.5\linespacing}{\normalfont\scshape\centering\S}}
\makeatother

\usepackage[abbrev,msc-links,backrefs]{amsrefs}
\usepackage{doi}

\renewcommand{\PrintDOI}[1]{\doi{#1}}

\usepackage{datetime}
\usepackage{lineno}
\newcommand*\patchAmsMathEnvironmentForLineno[1]{%
\expandafter\let\csname old#1\expandafter\endcsname\csname #1\endcsname
\expandafter\let\csname oldend#1\expandafter\endcsname\csname end#1\endcsname
\renewenvironment{#1}%
{\linenomath\csname old#1\endcsname}%
{\csname oldend#1\endcsname\endlinenomath}}%
\newcommand*\patchBothAmsMathEnvironmentsForLineno[1]{%
\patchAmsMathEnvironmentForLineno{#1}%
\patchAmsMathEnvironmentForLineno{#1*}}%
\AtBeginDocument{%
\patchBothAmsMathEnvironmentsForLineno{equation}%
\patchBothAmsMathEnvironmentsForLineno{align}%
\patchBothAmsMathEnvironmentsForLineno{flalign}%
\patchBothAmsMathEnvironmentsForLineno{alignat}%
\patchBothAmsMathEnvironmentsForLineno{gather}%
\patchBothAmsMathEnvironmentsForLineno{multline}%
}

\usepackage{geometry}
\newtheorem{theorem}{Theorem}
\newtheorem{lemma}[theorem]{Lemma}
\newtheorem{claim}[theorem]{Claim}

\newtheorem{corollary}[theorem]{Corollary}

\newtheorem{conjecture}[theorem]{Conjecture}
\newtheorem{proposition}[theorem]{Proposition}
\newtheorem*{question*}{Question}

\newcommand{\defi}[1]{%
  \emph{\color{red!60!black}#1}%
 }

\title{Biclique immersions in graphs with independence number 2}

\author{F. Botler}
\address[F. Botler]{
  Departamento de Ciência da Computação -- 
  Instituto de Matemática e Estatística -- 
  Universidade de São Paulo -- 
  Brasil
}

\author{A. Jiménez}
\address[A. Jiménez, D. A. Quiroz]{
  Instituto de Ingenier\'ia Matem\'atica \& CIMFAV --
  Universidad de Valpara\'iso -- 
  Chile
}

\author{C. N. Lintzmayer}
\address[C. N. Lintzmayer, M. Sambinelli]{
  Centro de Matemática, Computação e Cognição -- 
  Universidade Federal do ABC -- 
  Brazil
}

\author{A. Pastine}
\address[A. Pastine]{
  Instituto de Matemática Aplicada San Luis -- 
  CONICET \& Universidad Nacional de San Luis -- 
  Argentina
}

\author{D. A. Quiroz}

\author{M. Sambinelli}

\email{fbotler@ime.usp.br, andrea.jimenez@uv.cl, carla.negri@ufabc.edu.br, agpastine@unsl.edu.ar, daniel.quiroz@uv.cl, m.sambinelli@ufabc.edu.br}

\thanks{
  This research has been partially supported by Coordenação de Aperfeiçoamento
  de Pessoal de Nível Superior -- Brasil -- CAPES -- Finance Code 001, by
  MATHAMSUD MATH190013, and by FAPESP (Proc.~2019/13364-7).
 F.~Botler is supported by CNPq {(\small 304315/2022-2)} and
CAPES {(\small 88887.878880/2023-00)}
  A. Jiménez is partially supported by ANID/Fondecyt Regular 1220071, MATHAMSUD
  MATH210008, and ANID/PCI/REDES 190071.
  C. N. Lintzmayer is partially supported by CNPq (Proc.~312026/2021-8).
  A. Pastine  is partially supported by Universidad Nacional de San Luis,
  Argentina, grants PROICO 03-0918 and PROIPRO 03-1720, and by ANPCyT grants
  PICT-2020-SERIEA-04064 and PICT-2020-SERIEA-00549.
  D.A. Quiroz is partially supported by ANID/Fondecyt Iniciación en
  Investigación 11201251 and MATHAMSUD MATH210008.
  FAPESP is the Research Foundations S\~ao Paulo.  CNPq is the National Council for Scientific and
  Technological Development of Brazil.
}

\begin{document}

\maketitle

\begin{abstract}
  The analogue of Hadwiger's conjecture for the immersion relation states that
  every graph~$G$ contains an immersion of $K_{\chi(G)}$.
  For graphs with independence number 2, this is equivalent to stating that every
  such $n$-vertex graph contains an immersion of $K_{\lceil n/2 \rceil}$.
  We show that every $n$-vertex graph with independence number 2 contains every
  complete bipartite graph on $\lceil n/2 \rceil$ vertices as an immersion. 
\end{abstract}

\section{Introduction}
\label{sec:introduction}

A central problem in graph theory is guaranteeing dense substructures in graphs with a
given chromatic number.
Hadwiger's Conjecture~\cite{Hadwiger1943} is one of the most important examples of this pursuit,
stating that every loopless graph~$G$ contains the complete graph $K_{\chi(G)}$ as a minor 
(where \(\chi(G)\) is the chromatic number of \(G\)),
thus aiming to generalize the Four Color Theorem.
This difficult conjecture is known to hold whenever $\chi(G) \le 6$~\cite{Had6}, and it is open
for the remaining values. 
Thus, a natural approach is to study whether it holds whenever~$G$ is restricted to
a particular class of graphs.
A class of graphs that has received particular attention (and yet remains open) is that of 
graphs with independence number~2.
A recent survey of Seymour~\cite{Seymour2016} emphasizes the importance of this case, which
was first remarked by Mader (see~\cite{PlummerST2003}). 
Plummer, Stiebitz, and Toft~\cite{PlummerST2003} gave an equivalent formulation 
of Hadwiger's Conjecture for such graphs: 
every $n$-vertex graph with independence number~2 contains a minor
of $K_{\lceil n/2 \rceil}$.
Before that, in 1982, Duchet and Meyniel~\cite{DuchetM1982} had shown a result that implies that
every such graph contains a minor of $K_{\lceil n/3 \rceil}$.
Despite much work, see e.g.~\cite{Fox2010, KawaPT05, KawaSong2007, Wood2007}, it is still open
whether there is a constant $c > 1/3$ such that every graph with independence number~2 contains
a minor of $K_{\lceil cn \rceil}$.
Given the difficulty to obtain a complete minor on $\lceil n/2 \rceil$ vertices, Norin and
Seymour~\cite{NorinS2022} recently turned into finding dense minors on this amount of vertices.
They proved that every $n$-vertex graph with independence number~2 contains a (simple) minor of a graph~$H$ on $\lceil n/2 \rceil$ vertices and $0.98688 \cdot \binom{|V(H)|}{2} - o(n^2)$ edges.

The focus of this paper is a conjecture related to Hadwiger's, concerned with finding
graph immersions in graphs with a given chromatic number; this type of substructure is 
defined as follows.
To \defi{split off} a pair of adjacent edges $uv$, $vw$ (with \(u\), \(v\), \(w\) distinct) amounts to deleting those two edges 
and adding an edge joining \(u\) to \(v\) (keeping parallel edges if needed).
A graph~$G$ is said to contain an \defi{immersion} of another graph~$H$ if~$H$ can be obtained
from a subgraph of~$G$ by splitting off pairs of edges and deleting isolated vertices.
Notice then that if~$G$ contains~$H$ as a subdivision, it contains~$H$ as an immersion 
(and as a minor).
Immersions have received increased attention in recent years, see e.g.~\cite{DeVosMMS,DvorakWollan,Liu,LiuMuzi,MarxWollan,Wollan},
particularly since Robertson and Seymour~\cite{Robertson2010} proved that graphs are 
well-quasi-ordered by the immersion relation.
Much of this attention has been centered around the following conjecture of Lescure and Meyniel ~\cite{LescureM1989} (see also~\cite{Abu-KhzamL2003}), which is the immersion-analogue of Hadwiger's Conjecture.

\begin{conjecture}[Lescure and Meyniel \cite{LescureM1989}]
\label{Conj:AL}
  Every loopless graph~$G$ contains an immersion of~$K_{\chi(G)}$.
\end{conjecture}

The above conjecture holds whenever $\chi(G) \le 4$ because Hajós' subdivision conjecture holds in this case, actually giving a subdivision of $K_{\chi(G)}$~\cite{Dirac}. 
The cases where $\chi(G) \in \{5,6,7\}$ were verified independently by Lescure and
Meyniel~\cite{LescureM1989} and by DeVos, Kawarabayashi, Mohar, and Okamura~\cite{DeVosKMO2010}.
In general, a result of Gauthier, Le, and Wollan~\cite{GauthierLW2019}
guarantees that every graph~$G$ contains an immersion of a clique on 
$\lceil \frac{\chi(G) - 4}{3.54} \rceil$ vertices.
This result improves on theorems due to Dvo\v{r}ák and Yepremyan~\cite{dvovrak2018complete} 
and DeVos, Dvo\v{r}ák, Fox, McDonald, Mohar, and Scheide~\cite{DeVosDFMMS}.

The case of graphs with independence number~2 has also received attention in regard to Conjecture~\ref{Conj:AL}.
In particular, Vergara~\cite{Vergara2017} showed that, for such graphs, Conjecture~\ref{Conj:AL} is equivalent to the following conjecture.

\begin{conjecture}[Vergara \cite{Vergara2017}]
    Every $n$-vertex graph with independence number~2 contains an immersion of~$K_{\lceil n/2 \rceil}$.
\end{conjecture}

As evidence for her conjecture, Vergara proved that every $n$-vertex graph with independence
number~2 contains an immersion of~$K_{\lceil n/3 \rceil}$.
This was later improved by Gauthier \textit{et al.}~\cite{GauthierLW2019}, who showed that every 
such graph contains an immersion of $K_{2\lfloor n/5 \rfloor}$.
This last result was extended to graphs with arbitrary independence number~\cite{BustamanteQSZ2021}.
Additionally, Vergara’s Conjecture has been verified for graphs with small 
forbidden subgraphs~\cite{Quiroz2021}.
The main contribution of this paper is the following result, which states that graphs with 
independence number~2 contain an immersion of every complete bipartite graph on
$\lceil n/2 \rceil$ vertices.

\begin{theorem}
\label{theo:main}
    Let $G$ be an $n$-vertex graph with independence number~2, and $\ell \le \lceil n/2 \rceil - 1$
    be a positive integer.
    Then~$G$ contains an immersion of $K_{\ell, \lceil n/2 \rceil - \ell}$.
\end{theorem}

We remark that our proof of Theorem~\ref{theo:main} is self-contained. Using an argument due to Plummer \textit{et al.}~\cite{PlummerST2003} we show
(see Corollary~\ref{cor:equiv}) that Theorem~\ref{theo:main} implies the following.

\begin{corollary}
\label{cor:maincorollary}
    Let $G$ be a graph with independence number~2, and $1 \le \ell \le \chi(G) - 1$.
    Then~$G$ contains an immersion of $K_{\ell, \chi(G) - \ell}$.
\end{corollary}

This result leads us to make the following conjecture, which is a weakining of Conjecture \ref{Conj:AL}
and holds trivially when $\ell = 1$.

\begin{conjecture}
\label{conj:chromatic}
    If $1 \le \ell \le \chi(G) - 1$, then $G$ contains an immersion of $K_{\ell, \chi(G) - \ell}$.
\end{conjecture}

We denote by \defi{$K_{a, b, c}$} the graph that admits a partition into parts
of sizes $a$, $b$, and $c$ such that any pair of these parts induces a complete
bipartite graph.
In addition to Corollary~\ref{cor:maincorollary}, as evidence for
Conjecture~\ref{conj:chromatic}, we give a short proof of the following
strengthening of the case $\ell = 2$,
which is also implied by a special case of a result of Mader for subdivisions~\cite{mader1971existenz}.

\begin{proposition}
\label{prop:ellequals2}
    If $\chi(G) \geq 3$, then $G$ contains $K_{1, 1, \chi(G) - 2}$ as an immersion.
\end{proposition}

We note that Conjecture~\ref{conj:chromatic} has its parallel in the minor
order.
Woodall~\cite{Woodall2001} and, independently, Seymour (see~\cite{KostochkaP2010}),
proposed the following conjecture: every graph~$G$ with $\ell \le \chi(G) - 1$ 
contains a minor of $K_{\ell, \chi(G) - \ell}$.
In~\cite{Woodall2001}, Woodall showed that (the list-coloring strengthening of)
his conjecture holds whenever $\ell \le 2$. 
Kostochka and Prince~\cite{KostochkaP2010} showed that the case $\ell = 3$ holds
as long as $\chi(G) \ge 6503$. 
Kostochka~\cite{Kostochka2010} proved it for
every~$\ell$ as long as~$\chi(G)$ is very large in comparison to~$\ell$, 
and later~\cite{Kostochka2014} improved this so that~$\chi(G)$ could be
polynomial in~$\ell$, namely, whenever
$\chi(G) > 5 (200 \ell \log_2(200\ell))^3 + \ell$.
In fact, the results in~\cite{KostochkaP2010,Kostochka2010,Kostochka2014} obtain the full join
 $K_{\ell,\chi(G)-\ell}^*$, which is the graph obtained from the disjoint union of a \(K_\ell\)
 and an independent set on \(\chi(G) -\ell\) vertices by adding all of the possible edges between them.
This and the above-cited result of Norin and Seymour leads us to make the 
following conjecture, which can be shown to be (using the arguments of Corollary~\ref{cor:equiv}) a particular case of the conjecture of Woodall and Seymour.


\begin{conjecture}\label{conj:minor}
    Let $G$ be an $n$-vertex graph with independence number 2, and 
    $1 \le \ell \le \lceil n/2 \rceil - 1$.
    Then $G$ contains a minor of $K_{\ell, \lceil n/2 \rceil - \ell}$.
\end{conjecture}

Note that the result of Kostochka leaves open the balanced case, thus not implying Conjecture~\ref{conj:minor}. 
Yet, inspired by an earlier version of this paper Chen and Deng~\cite{chen2024seymour} showed that Conjecture~\ref{conj:minor} holds.

An extended abstract, containing a short sketch of the proof of Theorem~\ref{theo:main}, has appeared in~\cite{BoJiLiPaQuiSa23arxiv}.

The rest of the paper is organized as follows.
In Section~\ref{sec:preliminaries} we give some definitions and prove an
interesting lemma that is used to prove Theorem~\ref{theo:main}, in
Section~\ref{sec:main}.
In Section~\ref{sec:consequences} we prove Corollary~\ref{cor:maincorollary}
and Proposition~\ref{prop:ellequals2}.
Finally, in Section~\ref{sec:matching-immersions} we present a related result
together with an open problem.

\section{Preliminaries and notation}
\label{sec:preliminaries}

All the graphs considered in this work are finite and loopless.
For a graph~$G$, we denote by \defi{$V(G)$} and \defi{$E(G)$} its set of
vertices and edges, respectively. For every positive integer~$n$, let $\defi{[n]} = \{1, \dots , n\}$.
Let~$G$ be a graph.
For $v \in V(G)$ and $S \subseteq V(G)$, we define 
\defi{$E(v,S) = \{vu \in E(G) \colon u \in S\}$}.
For a set of vertices $W \subseteq V(G)$, 
we write \defi{$G[W]$} to denote the subgraph induced by~$W$.
For a set of edges $F \subseteq E(G)$, 
we write \defi{$G[F]$} to denote the subgraph induced by~$F$,
i.e., the subgraph of \(G\) whose edge set is \(F\),
and whose vertex set is the set of all end vertices of the edges in \(F\).
We also denote by \defi{$N_G(v)$} and \defi{$N_G[v]$} the neighborhood and the
closed neighborhood, respectively, of~$v$ in~$V(G)$.
The degree of a vertex~$v$ of a graph~$G$ is denoted by~\defi{$d_G(v)$}.
We simply write~$N(v)$, $N[v]$, and~$d(v)$ when~$G$ is clear from the context.
Moreover, when convenient, given a set $F \subseteq E(G)$ and a vertex $v \in
V(G)$, we write~$N_F(v)$, $N_F[v]$, and~$d_F(v)$ to denote, respectively,
$N_{G[F]}(v)$, $N_{G[F]}[v]$, and~$d_{G[F]}(v)$.

Given an edge~$e$ of a graph~$G$, the \defi{multiplicity of (the edge)~$e$ 
in~$G$}, denoted by~\defi{$m_G(e)$}, is the number of edges in~$G$ having the
same pair of vertices as the edge~$e$.
For sets $A,B \subseteq V(G)$, we say that a path~$P$ of~$G$ is an 
\defi{$A,B$-path} if~$P$ has an endvertex in~$A$ and the other in~$B$.
We denote by~\defi{$\alpha(G)$} the maximum size of an independent set of~$G$,
which is the \defi{independence number} of~$G$.

The \defi{complete graph} on~$n$ vertices is denoted by~\defi{$K_n$} and
the \defi{complete bipartite graph} with parts of size $a$ and $b$
is denoted by~\defi{$K_{a,b}$}. If~$A$ and~$B$ are disjoint sets, we let \defi{$K_{A,B}$} be the complete
bipartite graph with bipartition $(A,B)$.
Complete graphs and complete bipartite graphs are also called \defi{cliques} and \defi{bicliques}, respectively.


In a manner that is equivalent to the definition given in the introduction, given a graph \(H\), we
say that a graph~$G$ \defi{contains an immersion of~$H$} if there exists an
injection $f \colon V(H) \to V(G)$ and a collection of edge-disjoint paths 
in~$G$, one for each edge of~$H$, such that the path corresponding to 
an edge joining \(u\) and \(v\) has endvertices $f(u)$ and $f(v)$. 

In our proof, we make use of the following lemma, which we believe to be
interesting by itself, and that we could not find in the literature.
The lemma below may be rephrased in terms of (disjoint) systems of distinct
representatives.

\begin{figure}[h]
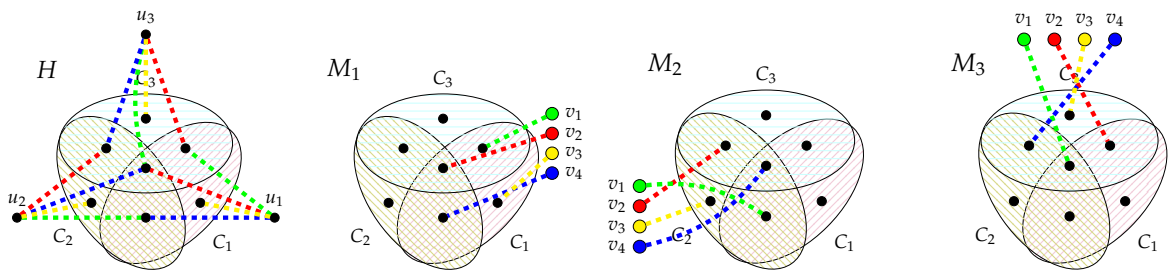

  \centering
  \resizebox{0.49\textwidth}{!}{
  \tikz{
  \node[ellipse,rotate=45,draw,minimum height=2cm, minimum width=3.5cm,label=below:{$C_1$},pattern=north east lines,pattern color=purple!20] at (.5,-1) (C1) {};
  \node[ellipse,draw,rotate=305,minimum height=2cm, minimum width=3.5cm,pattern=north west lines,pattern color=olive!30,label=below:{$C_2$}] at (-.5,-1) (C2) {};
  \node[ellipse,draw,label=above:{$C_3$},pattern=horizontal lines,pattern color=cyan!20,minimum height=2cm, minimum width=3.5cm] at (0,0) (C3) {};

  \node[circle,draw,fill,scale=.5] at (0,-.5) (c123) {};
  \node[circle,draw,fill,scale=.5] at (0,-1.5) (c12) {};
  \node[circle,draw,fill,scale=.5] at (.8,-.1) (c13) {};
  \node[circle,draw,fill,scale=.5] at (-.8,-.1) (c23) {};
  \node[circle,draw,fill,scale=.5] at (0,.5) (c3) {};
  \node[circle,draw,fill,scale=.5] at (1.1,-1.2) (c1) {};
  \node[circle,draw,fill,scale=.5] at (-1.1,-1.2) (c2) {};

  \node[circle,draw,fill,label=above:{$u_1$},scale=.5] at (2.6,-1.5) (u1) {};
  \node[circle,draw,fill,label=above:{$u_2$},scale=.5] at (-2.6,-1.5) (u2) {};
  \node[circle,draw,fill,label=above:{$u_3$},scale=.5] at (0,2.2) (u3) {};

  \draw (u1) [dashed,red,line width=2.5] to (c123);
  \draw (u1) [dashed,green,line width=2.5] to (c13);
  \draw (u1) [dashed,blue,line width=2.5] to (c12);
  \draw (u1) [dashed,yellow,line width=2.5] to (c1);

  \draw (u2) [dashed,blue,line width=2.5] to (c123);
  \draw (u2) [dashed,red,line width=2.5] to (c23);
  \draw (u2) [dashed,green,line width=2.5] to (c12);
  \draw (u2) [dashed,yellow,line width=2.5] to (c2);

  \draw (u3) [dashed,green,line width=2.5, bend right=15] to (c123);
  \draw (u3) [dashed,red,line width=2.5] to (c13);
  \draw (u3) [dashed,blue,line width=2.5] to (c23);
  \draw (u3) [dashed,yellow,line width=2.5] to (c3);

  \node[circle,scale=1.3] at (-2,1.5) (h) {$H$};

  \node[ellipse,rotate=45,draw,minimum height=2cm, minimum width=3.5cm,pattern=north east lines,pattern color=purple!20,label=below:{$C_1$}] at (6.5,-1) (C1) {};
  \node[ellipse,draw,rotate=305,minimum height=2cm, minimum width=3.5cm,pattern=north west lines,pattern color=olive!30,label=below:{$C_2$}] at (5.5,-1) (C2) {};
  \node[ellipse,draw,label=above:{$C_3$},pattern=horizontal lines,pattern color=cyan!20,minimum height=2cm, minimum width=3.5cm] at (6,0) (C3) {};

  \node[circle,draw,fill,scale=.5] at (6,-.5) (c123) {};
  \node[circle,draw,fill,scale=.5] at (6,-1.5) (c12) {};
  \node[circle,draw,fill,scale=.5] at (6.8,-.1) (c13) {};
  \node[circle,draw,fill,scale=.5] at (5.2,-.1) (c23) {};
  \node[circle,draw,fill,scale=.5] at (6,.5) (c3) {};
  \node[circle,draw,fill,scale=.5] at (7.1,-1.2) (c1) {};
  \node[circle,draw,fill,scale=.5] at (4.9,-1.2) (c2) {};

  \node[circle,draw,fill=green,label=right:{$v_1$},scale=.7] at (8.2,.6) (a1) {};
  \node[circle,draw,fill=red,label=right:{$v_2$},scale=.7] at (8.2,.2) (a2) {};
  \node[circle,draw,fill=yellow,label=right:{$v_3$},scale=.7] at (8.2,-.2) (a3) {};
  \node[circle,draw,fill=blue,label=right:{$v_4$},scale=.7] at (8.2,-.6) (a4) {};

  \node[circle,scale=1.3] at (4,1.5) (m1) {$M_1$};

  \draw (a2) [dashed,red,line width=2.5] to (c123);
  \draw (a1) [dashed,green,line width=2.5] to (c13);
  \draw (a4) [dashed,blue,line width=2.5] to (c12);
  \draw (a3) [dashed,yellow,line width=2.5] to (c1);
  }}
  \resizebox{0.48\textwidth}{!}{
  \tikz{
  \node[ellipse,rotate=45,draw,minimum height=2cm, minimum width=3.5cm,pattern=north east lines,pattern color=purple!20,label=below:{$C_1$}] at (.5,-1) (C1) {};
  \node[ellipse,draw,rotate=305,minimum height=2cm, minimum width=3.5cm,pattern=north west lines,pattern color=olive!30,label=below:{$C_2$}] at (-.5,-1) (C2) {};
  \node[ellipse,draw,label=above:{$C_3$},pattern=horizontal lines,pattern color=cyan!20,minimum height=2cm, minimum width=3.5cm] at (0,0) (C3) {};

  \node[circle,draw,fill,scale=.5] at (0,-.5) (c123) {};
  \node[circle,draw,fill,scale=.5] at (0,-1.5) (c12) {};
  \node[circle,draw,fill,scale=.5] at (.8,-.1) (c13) {};
  \node[circle,draw,fill,scale=.5] at (-.8,-.1) (c23) {};
  \node[circle,draw,fill,scale=.5] at (0,.5) (c3) {};
  \node[circle,draw,fill,scale=.5] at (1.1,-1.2) (c1) {};
  \node[circle,draw,fill,scale=.5] at (-1.1,-1.2) (c2) {};

  \node[circle,draw,fill=green,label=left:{$v_1$},scale=.7] at (-2.5,-.9) (a1) {};
  \node[circle,draw,fill=red,label=left:{$v_2$},scale=.7] at (-2.5,-1.3) (a2) {};
  \node[circle,draw,fill=yellow,label=left:{$v_3$},scale=.7] at (-2.5,-1.7) (a3) {};
  \node[circle,draw,fill=blue,label=left:{$v_4$},scale=.7] at (-2.5,-2.1) (a4) {};

  \node[circle,scale=1.3] at (-2,1.5) (m2) {$M_2$};

  \draw (a2) [dashed,red,line width=2.5] to (c23);
  \draw (a1) [dashed,green,line width=2.5, bend left = 20] to (c12);
  \draw (a4) [dashed,blue,line width=2.5, bend right=20] to (c123);
  \draw (a3) [dashed,yellow,line width=2.5] to (c2);

  \node[ellipse,rotate=45,draw,minimum height=2cm,pattern=north east lines,pattern color=purple!20, minimum width=3.5cm,label=below:{$C_1$}] at (6.5,-1) (C1) {};
  \node[ellipse,draw,rotate=305,minimum height=2cm,pattern=north west lines,pattern color=olive!30, minimum width=3.5cm,label=below:{$C_2$}] at (5.5,-1) (C2) {};
  \node[ellipse,draw,label=above:{$C_3$},pattern=horizontal lines,pattern color=cyan!20,minimum height=2cm, minimum width=3.5cm] at (6,0) (C3) {};

  \node[circle,draw,fill,scale=.5] at (6,-.5) (c123) {};
  \node[circle,draw,fill,scale=.5] at (6,-1.5) (c12) {};
  \node[circle,draw,fill,scale=.5] at (6.8,-.1) (c13) {};
  \node[circle,draw,fill,scale=.5] at (5.2,-.1) (c23) {};
  \node[circle,draw,fill,scale=.5] at (6,.5) (c3) {};
  \node[circle,draw,fill,scale=.5] at (7.1,-1.2) (c1) {};
  \node[circle,draw,fill,scale=.5] at (4.9,-1.2) (c2) {};

  \node[circle,draw,fill=green,label=above:{$v_1$},scale=.7] at (5.1,2) (a1) {};
  \node[circle,draw,fill=red,label=above:{$v_2$},scale=.7] at (5.7,2) (a2) {};
  \node[circle,draw,fill=yellow,label=above:{$v_3$},scale=.7] at (6.3,2) (a3) {};
  \node[circle,draw,fill=blue,label=above:{$v_4$},scale=.7] at (6.9,2) (a4) {};

  \node[circle,scale=1.3] at (4,1.5) (m1) {$M_3$};

  \draw (a2) [dashed,red,line width=2.5] to (c13);
  \draw (a1) [dashed,green,line width=2.5] to (c123);
  \draw (a4) [dashed,blue,line width=2.5] to (c23);
  \draw (a3) [dashed,yellow,line width=2.5] to (c3);
  }}
  
  \caption{An example of Lemma~\ref{lem:matching_lemma}.}
  \label{fig:matching lemma}
\end{figure}

\begin{lemma}
\label{lem:matching_lemma}
  Let $j \leq k$ be two positive integers, and let $C_1, \ldots, C_j \subseteq [n]$
  be sets of size~$k$.
  Let~$A$ be a set of size~$k$ disjoint from~$[n]$.
  Then there are disjoint matchings $M_1, \ldots, M_j$ in $K_{A,[n]}$ such 
  that~$M_i$ matches~$A$ with~$C_i$ for every $i \in [j]$.
\end{lemma}
\begin{proof}
  Let \(H\) be an auxiliary bipartite graph such that 
  $V(H) = \{u_1, \ldots, u_j\} \cup [n]$ and $E(H) = \{u_ix \colon x \in C_i, i \in [j]\}$.
  Since $|C_i| = k$, we have $d(u_i) = k$.
  Moreover, the degree of a vertex $x \in [n]$ is precisely the number of $C_i$'s
  that contain~$x$, which is at most $j \leq k$.
  Therefore, $H$ is a bipartite graph of maximum degree precisely~$k$.
  Then, by K\"onig's line coloring theorem, we can find a proper \(k\)-coloring of the edges of~$H$, say $D_1, \ldots, D_k \subseteq E(H)$.

  Let $A = \{v_1, \ldots, v_k\}$ be disjoint from~$[n]$.
  Now we construct the matchings of $K_{A,[n]}$.
  For each $i \in \{1, \ldots, k\}$, let $u_i w^i_1, \ldots, u_i w^i_k$ be the
  edges of~$H$ incident to~$u_i$, where $u_i w^i_j \in D_j$ (for an example, see
  Figure~\ref{fig:matching lemma}).
  Note that, by the definition of~$H$, we have $w^i_1, \ldots, w^i_k \in C_i$.
  Then we create the matching $M_i = \{v_1 w^i_1, \ldots, v_k w^i_k\}$, which
  matches~$A$ with~$C_i$.

  We claim that $M_i \cap M_j = \emptyset$ for every $i \neq j$.
  Suppose, without loss of generality, that $v_1 w \in M_i \cap M_j$.
  By the construction of~$M_i$, we have that $u_i w \in D_1$.
  Analogously, $u_j w \in D_1$.
  Then, since~$D_1$ is a matching, $u_i = u_j$, which means $i = j$.
  Thus, $M_i \cap M_j = \emptyset$ when $i\neq j$.
\end{proof}

\section{Proof of Theorem~\ref{theo:main}}
\label{sec:main}

While Theorem~\ref{theo:main} holds even when parallel edges are allowed, we will only consider simple graphs in this proof as this suffices and makes the writing simpler.

Indeed we can consider graphs with independence number at most 2. The proof follows by induction on $n + \ell$.
Let $G$ be an $n$-vertex graph with $\alpha(G) \le 2$ and let 
$\ell \leq \lceil n/2 \rceil - 1$ be a positive integer.
It is not hard to check that the result holds when $n\le 4$, and thus we may assume that $n\ge 5$.
Note also that it suffices to prove the statement in the case~$G$ is
\emph{edge-critical}, i.e., that the removal of any edge of~$G$ increases its
independence number.
Now, if $n \leq 4\ell - 2$, 
then $\lceil n/2 \rceil - \ell \leq 2\ell - 1 - \ell < \ell$.
Thus, by induction, there is an immersion of 
$K_{\ell', \lceil n/2 \rceil - \ell'}$ in~$G$, 
where $\ell' = \lceil n/2 \rceil - \ell$.
But this is the desired immersion because 
$K_{\lceil n/2 \rceil - \ell, \lceil n/2 \rceil - \lceil n/2 \rceil + \ell}=K_{\lceil n/2 \rceil - \ell, \ell}$.
Thus, from now on, we may assume that 
\begin{equation}
\label{eq:size}
    n \geq 4\ell - 1.
\end{equation}

Our proof is divided into four parts, and the following part is the only step of the proof
in which we directly use induction.

\subsection{Non-adjacent vertices with at least $\ell-1$ common neighbors}

\begin{figure}[h]
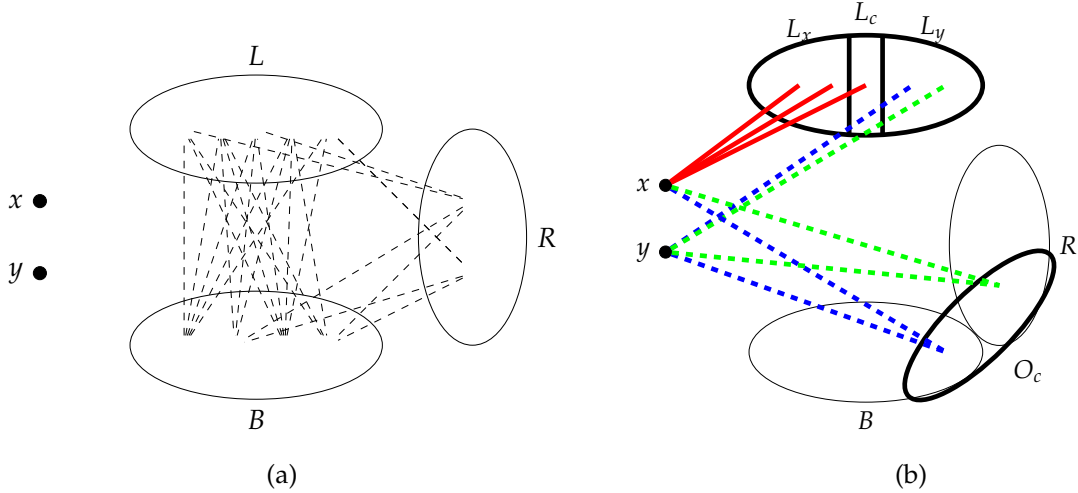

  \centering
  \begin{subfigure}[b]{0.48\textwidth}
    \resizebox{\textwidth}{!}{
    \tikz{
      \node[ellipse,draw,label=above:{$L$},minimum height=1.5cm, minimum width=3.5cm] at (3,1.5) (L) {};
      \node[ellipse,draw,label=below:{$B$},minimum height=1.5cm, minimum width=3.5cm] at (3,-1.5) (N-L) {};
      \node[ellipse,draw,label=right:{$R$},minimum height=3cm, minimum width=1.5cm] at (6,0) (rest) {};
  
      \node[circle,draw,fill,label=left:{$x$},scale=.5] at (0,0.5) (x) {};
      \node[circle,draw,fill,label=left:{$y$},scale=.5] at (0,-0.5) (y) {};
  
      \node at (2,1.5) (l1) {};
      \node at (2.5,1.5) (l2) {};
      \node at (3,1.5) (l3) {};
      \node at (3.5,1.5) (l4) {};
      \node at (4,1.5) (l5) {};
      
      \node at (2,-1.5) (b1) {};
      \node at (2.7,-1.5) (b2) {};
      \node at (3.4,-1.5) (b3) {};
      \node at (4,-1.5) (b4) {};
      
      \node at (6,0.5) (r1) {};
      \node at (6,-0.5) (r2) {};
  
      \foreach \l in {l1,l2,l3,l4,l5}{
          \foreach \b in {b1,b3}{
              \draw (\b) [dashed] to (\l);
            }
        }
  
      \draw (l3) [dashed] to (r1);
      \draw (r1) [dashed] to (b4);
      \draw (l5) [dashed] to (r2);
      \draw (r2) [dashed] to (b4);
      \draw (l1) [dashed] to (b4);
      \draw (l2) [dashed] to (b4);
      \draw (l4) [dashed] to (b4);
      
      \draw (l1) [dashed] to (r1);
      \draw (r1) [dashed] to (b2);
      \draw (l5) [dashed] to (r2);
      \draw (r2) [dashed] to (b2);
      \draw (l3) [dashed] to (b2);
      \draw (l2) [dashed] to (b2);
      \draw (l4) [dashed] to (b2);
      
    }
    }
    \caption{}
    \label{fig:induction_immersion}
  \end{subfigure}
  \hfill
  \begin{subfigure}[b]{0.48\textwidth}
    \resizebox{0.8\textwidth}{!}{
    \tikz{
    \node[ellipse,draw,label=above:{$L_c$},minimum height=1.5cm, minimum width=3.5cm,line width=2] at (3,2) (Lo) {};
    \node[draw=none] at (2,2.83) (Lx) {$L_x$};
    \node[draw=none] at (4,2.83) (Ly) {$L_y$};
    \draw (2.75,1.25) [line width=2] to (2.75,2.75);
    \draw (3.25,1.25) [line width=2] to (3.25,2.75);

    \node[ellipse,label=below:{$O_c$},rotate=45,draw,minimum height=1cm, minimum width=3cm,line width=2] at (4.7,-1.6) (Oc) {};

    \node[ellipse,draw,label=below:{$B$},minimum height=1.5cm, minimum width=3.5cm] at (3,-2) (N-L) {};
    \node[ellipse,draw,label=right:{$R$},minimum height=3cm, minimum width=1.5cm] at (5,-0.4) (rest) {};

    \node[circle,draw,fill,label=left:{$x$},scale=.5] at (0,0.5) (x) {};
    \node[circle,draw,fill,label=left:{$y$},scale=.5] at (0,-0.5) (y) {};

    \draw (x) [line width=2,color=red] to (2,2);
    \draw (x) [line width=2,color=red] to (2.5,2);
    \draw (x) [line width=2,color=red] to (3,2);
    \draw (y) [line width=2,dashed,color=blue] to (3.7,2);
    \draw (y) [line width=2,dashed,color=green] to (4.2,2);

    \draw (x) [line width=2,dashed,color=blue] to (4.2,-2);
    \draw (y) [line width=2,dashed,color=blue] to (4.2,-2);

    \draw (x) [line width=2,dashed,color=green] to (5,-1);
    \draw (y) [line width=2,dashed,color=green] to (5,-1);
    }}
    \caption{}
    \label{fig:paths_using_Oc}
  \end{subfigure}
  \caption{\eqref{fig:induction_immersion} Partition of~$V(G)$ with respect to~$H'$;
  \eqref{fig:paths_using_Oc} Paths from~$x$ to~$L_y$ using~$O_c$ and~$y$ depicted in dashed lines.}
\end{figure}

Our first task is to show that we can assume that 
\begin{equation}
\label{eq:common}
  |N(u) \cap N(v)| \leq \ell - 2 \textrm{ for every pair of non-adjacent vertices } u,v,
\end{equation}
as, otherwise, we have the desired immersion.
For that, suppose that~$G$ contains two non-adjacent vertices, say~$x$
and~$y$, with at least $\ell-1$ common neighbors, and let $G' = G - x - y$.
If $\ell \leq \lceil n/2 \rceil - 2 = \lceil (n-2)/2 \rceil - 1$ the induction
hypothesis guarantees that~$G'$ contains an immersion of 
$K_{\ell, \lceil (n-2)/2 \rceil - \ell}$, which we call~$H'$.
Otherwise, if we have $\ell = \lceil n/2 \rceil - 1$ we let~$H'$ be an 
arbitrary set of~$\ell$ vertices.
Let~$L$ and~$B$ be the parts (i.e., the sets of vertices of \(G\) on which the parts of $K_{\ell, \lceil (n-2)/2 \rceil - \ell}$ were mapped) of~$H'$ having size~$\ell$ and 
$\lceil n/2 \rceil - 1 - \ell$, respectively, and 
let~$R = V(G') \setminus (L \cup B)$ (see Figure~\ref{fig:induction_immersion}).
As~$\alpha(G) = 2$, every vertex in~$G'$ is either adjacent to~$x$ or to~$y$
in~$G$.
This is true in particular for the vertices in~$L$.
In what follows, we add either~$x$ or~$y$ to~$B$, in order to obtain the
desired immersion of $K_{\ell, \lceil n/2 \rceil - \ell}$.
This is immediate if~$x$ or~$y$ is adjacent to every vertex in~$L$.
Thus we may assume that $|E(y,L)|, |E(x,L)| < \ell$.
Now, consider the following sets.
\begin{itemize}
  \item[$\triangleright$] $L_x$, the set of vertices in $L$ adjacent to $x$ but not adjacent to $y$;
  \item[$\triangleright$] $L_y$, the set of vertices in $L$ adjacent to $y$ but not adjacent to $x$;
  \item[$\triangleright$] $L_c$, the set of vertices in $L$ adjacent to both $x$ and $y$; and
  \item[$\triangleright$] $O_c$, the set of vertices adjacent to both $x$ and $y$ that are not in $L$.
\end{itemize}

As~$x$ is adjacent to every vertex in $L_x \cup L_c$, it is enough to find
(edge-disjoint) paths from~$x$ to~$L_y$ without using edges of~$H'$.
Notice that $|L_y| + |L_c| = |E(y,L)| \leq \ell - 1$ and that, by hypothesis,
we have $|L_c| + |O_c| \geq \ell - 1$.
Thus $|O_c| \geq |L_y|$.
Let $O_c = \{o_1, o_2, \ldots, o_{|O_c|}\}$ and 
$L_y = \{\ell_1, \ell_2, \ldots, \ell_{|L_y|}\}$.
For $1 \leq i \leq |L_y|$, we take the path $x o_i y \ell_i$ 
(see Figure~\ref{fig:paths_using_Oc}).
These paths are clearly pairwise edge-disjoint and do not use edges of~$H'$
because the edges $x o_i$, $o_i y$, and $y \ell_i$ are not in~$G'$.
This concludes the proof of~\eqref{eq:common}.

\subsection{Direct consequences of edge-criticality}

Recall that~$G$ is edge-critical, meaning that the removal of any edge
$uv \in E(G)$ creates an independent set of size~3.
Hence, for such an edge there is a vertex~$w$ that is not adjacent to both
$u$ and~$v$.
We formalize this argument in the following claim.

\begin{claim}
\label{claim:edge-criticality}
  For any $u, v \in V(G)$, we have $uv \in E(G)$ if and only if 
  $N[u] \cup N[v] \neq V(G)$.
\end{claim}

For the rest of the proof, we fix two non-adjacent vertices~$x$ and~$y$, and
partition the vertices of~$G$ into the following three sets:
\begin{itemize}
  \item[$\triangleright$] $C = N(x) \cap N(y)$, the set of common neighbors of $x$ and $y$;
  \item[$\triangleright$] $X = \overline{N[y]}$, the set of non-neighbors of $y$ excluding $y$, which contains $x$; and
  \item[$\triangleright$] $Y = \overline{N[x]}$, the set of non-neighbors of $x$ excluding $x$, which contains $y$.
\end{itemize}
Note that by~\eqref{eq:common}, we have $|C| \leq \ell - 2$. 
Moreover, both~$X$ and~$Y$ induce complete subgraphs of~$G$, otherwise we could find an independent set of size~3.
This implies that \(|X|, |Y| < \lceil n/2\rceil\),
otherwise we could find $K_{\ell, \lceil (n-2)/2 \rceil - \ell}$ as a subgraph,
and hence \(C \neq\emptyset\).
Moreover, the edge-criticality of~$G$ yields the following claim.

\begin{claim}
\label{claim:no_v_in_C_sees_all_X_Y}
  For every vertex $a\in C$, we have $X,Y \nsubseteq N(a)$.
\end{claim}
\begin{proof}
  First, note that $N[x] = X \cup C$.
  Thus, since \(a\in N(x)\), if \(Y\subseteq N(a)\),
  then~$xa$ is an edge of~$G$ for which $N[x] \cup N[a] = V(G)$, 
  a contradiction to Claim~\ref{claim:edge-criticality}. 
  Analogously, we have $X \nsubseteq N(a)$.
\end{proof}

\begin{figure}[h]
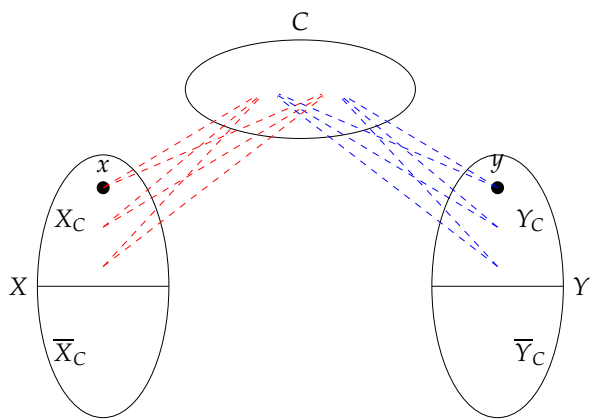

  \centering
  \resizebox{0.5\textwidth}{!}{
  \tikz{
  \node[ellipse,draw,label=above:{$C$},minimum height=1.5cm, minimum width=3.5cm] at (3,3) (C) {};
  \node[ellipse,draw,label=left:{$X$},minimum height=4cm, minimum width=2cm] at (0,0) (X) {};
  \node[ellipse,draw,label=right:{$Y$},minimum height=4cm, minimum width=2cm] at (6,0) (Y) {};
  \node[circle] at (2.5,3) (a) {};
  \node[circle] at (3.5,3) (b) {};
  \node[circle] at (-.5,1) (Ax) {$X_C$};
  \node[circle] at (6.5,1) (Ay) {$Y_C$};
  \node[circle] at (-.5,-1) (Bx) {$\overline{X}_C$};
  \node[circle] at (6.5,-1) (By) {$\overline{Y}_C$};
  \node[circle,draw,fill,label=above:{$x$},scale=.5] at (0,1.5) (x) {};
  \node[circle,draw,fill,label=above:{$y$},scale=.5] at (6,1.5) (y) {};

  \draw (-1,0) to (1,0);
  \draw (5,0) to (7,0);
  \foreach \i in {0,1,2}{
      \draw (0,1.5-.6*\i) [dashed, red] to (a);
      \draw (6,1.5-.6*\i) [dashed, blue] to (a);
      \draw (0,1.5-.6*\i) [dashed, red] to (b);
      \draw (6,1.5-.6*\i) [dashed, blue] to (b);
    }
  }}
  \caption{Partition of~$G$ according to non-adjacent vertices~$x$ and~$y$ and their common neighborhood~$C$.}
  \label{fig:partition_XXC_YYC}
\end{figure}

\subsection{Key vertex sets and their sizes}

We partition~$X$ and~$Y$ in different ways.
Let $X_C \subseteq X$ (resp.\ $Y_C \subseteq Y$) be the set containing vertices 
$v \in X$ (resp.\ $v \in Y$) for which $C \subset N(v)$, and put 
$\overline{X}_C = X \setminus X_C$ (resp.\ $\overline{Y}_C = Y \setminus Y_C$)
(see Figure~\ref{fig:partition_XXC_YYC}).
Similarly, given a vertex $a$ in $C$, we denote by~$X_a$ (resp.\ $Y_a$) the set of
vertices in~$X$ (resp.\ in~$Y$) that are adjacent to~$a$, and put
$\overline{X}_a = X \setminus X_a$ and $\overline{Y}_a = Y \setminus Y_a$.
Notice that if $v \in \overline{X}_a$ and $w \in \overline{Y}_a$, then~$v$
and~$w$ must be adjacent, as the independence number of~$G$ is~2.
Thus we get $K_{\overline{X}_a, \overline{Y}_a}$ as a subgraph of~$G$.
Note that $X_C \subseteq X_a$ and $\overline{X}_C \supseteq \overline{X}_a$
(resp.\ $Y_C \subseteq Y_a$ and $\overline{Y}_C \supseteq \overline{Y}_a$) 
for every $a \in C$.
Indeed, we have $X_C = \bigcap_{a \in C} X_a$ 
(resp.\ $Y_C = \bigcap_{a \in C} Y_a$)
and $\overline{X}_C = \bigcup_{a \in C} \overline{X}_a$ 
(resp.\ $\overline{Y}_C = \bigcup_{a \in C} \overline{Y}_a$).
The following claims give bounds on the sizes of some of these sets.
This control is the key to build the desired immersion.

\begin{claim}
  \label{claim:bounds_Xa_Ya}
  Given $a \in C$, we have that $|X_C| \leq |X_a| \leq \ell - 2$ and 
  $|Y_C| \leq |Y_a| \leq \ell - 2$.
  Furthermore, we have $|\overline{X}_a| \geq \lceil n/2 \rceil - |Y| + 3$ and
  $|\overline{Y}_a| \geq \lceil n/2 \rceil - |X| + 3$.
\end{claim}
\begin{proof}
  We prove the bounds for~$X_a$ and $\overline{X}_a$, as the remaining bounds
  follow analogously.
  By Claim~\ref{claim:no_v_in_C_sees_all_X_Y},
  there is a vertex $w \in \overline{X}_C$ such that $wa \notin E(G)$.
  By the definition of~$X_a$ we have $X_C \subseteq X_a \subseteq N(a)$, which
  gives us $|X_C| \leq |X_a|$.
  Moreover, since~$X$ induces a complete graph, we have $X_a \subseteq N(w)$,
  and hence $X_a \subseteq N(a) \cap N(w)$.
  Therefore, \eqref{eq:common} implies $|X_a| \leq \ell - 2$.

  Now, we use \eqref{eq:size} in the following to obtain the claim:
  \begin{align*}
        |Y| + |\overline{X}_a| &= n - |C| - |X_a| \\
                               &\geq n - \ell+2 - \ell+2 \\
                               &= \lceil n/2 \rceil + \lfloor n/2 \rfloor - 2\ell + 4 \\
                               &\geq \lceil n/2 \rceil + \lfloor (4\ell-1)/2 \rfloor - 2\ell + 4 \\
                               &= \lceil n/2 \rceil + \lfloor 2\ell - (1/2) \rfloor - 2\ell + 4 \\
                               &= \lceil n/2 \rceil + 2\ell - 1 - 2\ell + 4 \\
                               &= \lceil n/2 \rceil + 3 \, . \qedhere
  \end{align*}
\end{proof}

The claim above is used to prove the following three claims.
In fact, the ``+3'' term in its bounds is almost necessary for the proof of the
next claim.

\begin{claim}
\label{claim:y-y(a) is big}
  For every $a \in C$, we have either 
  $|X_a| \geq 2 (\ell - |\overline{Y}_a|)$ or 
  $|Y_a| \geq 2 (\ell - |\overline{X}_a|)$.
\end{claim}
\begin{proof}
  For a contradiction, assume 
  $|X_a| < 2 (\ell -|\overline{Y}_a|)$ and
  $|Y_a| < 2 (\ell -|\overline{X}_a|)$.
  Then we have 
  \begin{align*}
        |X| + |Y|
          &= |X| - |\overline{X}_a| + |\overline{X}_a| + |Y| - |\overline{Y}_a| + |\overline{Y}_a| \\
          &< 2 (\ell - |\overline{Y}_a|) + |\overline{X}_a| + 2 (\ell -|\overline{X}_a|) + |\overline{Y}_a| \\
          &= 4\ell - |\overline{Y}_a| - |\overline{X}_a| \, .
  \end{align*} 
  Also, by Claim~\ref{claim:bounds_Xa_Ya}, we have
  $|\overline{Y}_a| \geq \lceil n/2 \rceil - |X| + 3$ and 
  $|\overline{X}_a| \geq \lceil n/2 \rceil - |Y| + 3$, which implies 
  \begin{align*}
      4\ell & > |X| + |Y| + |\overline{Y}_a| + |\overline{X}_a| \\
          &\geq |X| + |Y| + \lceil n/2 \rceil - |X| + 3 + \lceil n/2 \rceil - |Y| + 3 \\
          &\geq n + 6\,.
  \end{align*} 
  So we have $n < 4\ell - 6$, which contradicts~\eqref{eq:size}.
\end{proof}

\begin{claim}
\label{claim:neighborhoods_Xprime_Yprime}
  For every $v \in \overline{X}_C$ (resp.\ $w \in \overline{Y}_C$), we have 
  $|N(v) \cap \overline{Y}_C| > \lceil n/2 \rceil - |X|$
  (resp.\ $|N(w) \cap \overline{X}_C| > \lceil n/2 \rceil - |Y|$).
\end{claim}
\begin{proof}
  We prove the bound for $|N(v) \cap \overline{Y}_C|$, as the bound on
  $|N(w) \cap \overline{X}_C|$ follows analogously.
  By the definition of $\overline{X}_C$, there is a vertex $a \in C$ such that $va \notin E(G)$.
  Note that $\overline{Y}_a \subseteq N(v)$, otherwise~$G$ would contain an
  independent set of size~3.
  Since $\overline{Y}_a \subseteq \overline{Y}_C$, we have $\overline{Y}_a
  \subseteq N(v) \cap \overline{Y}_C$, and the desired bound follows from
  Claim~\ref{claim:bounds_Xa_Ya}.
\end{proof}

Finally, we need the following relation between the sizes of~$\overline{X}_C$
and~$X_C$ (resp.\ $\overline{Y}_C$ and~$Y_C$).

\begin{claim}
\label{claim:if y' is small, then y^c is big}
  For every $a \in C$ we have the following:
  \begin{enumerate}[i)]
    \item If $|\overline{X}_C| < \ell$, then
    $|X_C| > \lceil n/2 \rceil - \ell - |\overline{Y}_a| \geq \ell - |\overline{Y}_a|$; and 
    \label{claim:if y' is small, then y^c is big:1}
    
    \item If $|\overline{Y}_C| < \ell$, then
    $|Y_C| > \lceil n/2 \rceil - \ell - |\overline{X}_a| \geq \ell - |\overline{X}_a| $.
    \label{claim:if y' is small, then y^c is big:2}
  \end{enumerate}
\end{claim}
\begin{proof}
  We only prove~\eqref{claim:if y' is small, then y^c is big:1},
  as~\eqref{claim:if y' is small, then y^c is big:2} follows analogously.
  By Claim~\ref{claim:bounds_Xa_Ya}, if $|\overline{X}_C| < \ell$, then
  $$\lceil n/2 \rceil + 3 \leq |\overline{Y}_a| + |X| = |\overline{Y}_a| +
    |\overline{X}_C| + |X_C| < |\overline{Y}_a| + \ell + |X_C| \, . $$

  By \eqref{eq:size} we have $\lceil n/2 \rceil \geq 2\ell$, which implies
  \[ |X_C| > \lceil n/2 \rceil - \ell - |\overline{Y}_a| \geq 2 \ell - \ell -
     |\overline{Y}_a| = \ell - |\overline{Y}_a| \, .\qedhere\]
\end{proof}

\subsection{Constructing the immersion}

The rest of the proof is divided into two cases which depend on the sizes 
of~$\overline{X}_C$ and~$\overline{Y}_C$.
The sets that form the bipartition of the immersion depend on which case we are
dealing with.
The construction requires more care in Case 2, where one of
$\overline{X}_C, \overline{Y}_C$ is large.
In that case, we apply Lemma~\ref{lem:matching_lemma} and obtain a set of paths
that are not necessarily edge-disjoint.
Nevertheless, the intersections are relatively few and with a combination of
different techniques, we are able to fix them and obtain the desired immersion.
For the rest of the proof, we fix $a \in C.$

\begin{figure}[h]
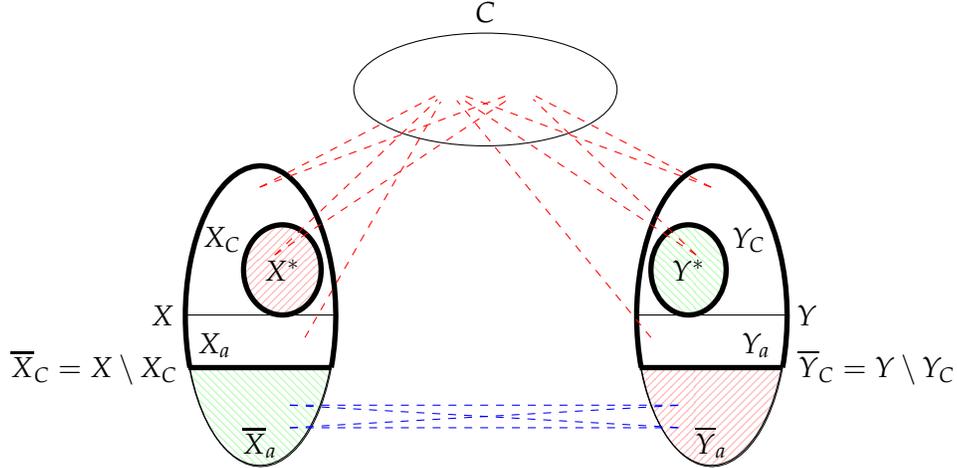

  \centering
  \tikz{
  \draw[pattern=north east lines,pattern color=red!30] (6.95,-0.7) arc(-20:-160:1cm and 2cm);
  \draw[pattern=north west lines,pattern color=green!30] (0.95,-0.7) arc(-20:-160:1cm and 2cm);

  \node[ellipse,draw,label=above:{$C$},minimum height=1.5cm, minimum width=3.5cm] at (3,3) (C) {};
  \node[ellipse,draw,label=left:{$X$},minimum height=4cm, minimum width=2cm] at (0,0) (X) {};
  \node[ellipse,draw,label=right:{$Y$},minimum height=4cm, minimum width=2cm] at (6,0) (Y) {};

  \node[ellipse,pattern=north east lines,pattern color=red!30,draw,minimum height=1.2cm, minimum width=1cm,line width=2pt] at (0.3,0.6) (Xs) {$X^*$};
  \node[ellipse,pattern=north west lines,pattern color=green!30,draw,minimum height=1.2cm, minimum width=1cm,line width=2pt] at (5.7,0.6) (Ys) {$Y^*$};

  \draw[line width=2pt] (0.95,-0.7) arc(-20:200:1cm and 2cm);
  \draw[line width=2pt] (-0.95,-0.7) to (0.95,-0.7);
  \node[circle] at (-0.6, -0.4) (Xa) {$X_a$};

  \draw[line width=2pt] (6.95,-0.7) arc(-20:200:1cm and 2cm);
  \draw[line width=2pt] (5.05,-0.7) to (6.95,-0.7);
  \node[circle] at (6.6, -0.4) (Ya) {$Y_a$};

  \node[circle] at (2.5,3) (a) {};
  \node[circle] at (3.5,3) (b) {};
  \node[circle] at (-.5,1) (XC) {$X_C$};
  \node[circle] at (6.5,1) (YC) {$Y_C$};
  \node[circle] at (0,-1.7) (Xab) {$\overline{X}_a$};
  \node[circle] at (6,-1.7) (Ybb) {$\overline{Y}_a$};

  \node at (-2.2,-0.7) (Xp) {$\overline{X}_C = X \setminus X_C$};
  \node at (8.2,-0.7) (Yp) {$\overline{Y}_C = Y \setminus Y_C$};

  \draw (-1,0) to (1,0);
  \draw (5,0) to (7,0);

  \draw (0,1.7) [dashed, red] to (a);
  \draw (0,1.7) [dashed, red] to (b);
  \draw (0.2,.8) [dashed, red] to (a);
  \draw (0.2,.8) [dashed, red] to (b);
  \draw (6,1.7) [dashed, red] to (a);
  \draw (6,1.7) [dashed, red] to (b);
  \draw (5.8,.8) [dashed, red] to (a);
  \draw (5.8,.8) [dashed, red] to (b);
  \draw (5.2,-.3) [dashed, red] to (a);
  \draw (0.6,-.3) [dashed, red] to (a);
  \draw (0.4,-1.5) [dashed, blue] to (5.6,-1.5);
  \draw (0.4,-1.2) [dashed, blue] to (5.6,-1.2);
  \draw (0.4,-1.5) [dashed, blue] to (5.6,-1.2);
  \draw (0.4,-1.2) [dashed, blue] to (5.6,-1.5);
  }
  \caption{Immersion for Case 1, when $|\overline{X}_C|, |\overline{Y}_C| < \ell$: 
  between $Y^* \cup \overline{X}_a$ and $X^* \cup \overline{Y}_a$.}
  \label{fig:case1}
\end{figure}

\smallskip
\noindent\textbf{Case 1}, $|\overline{X}_C|, |\overline{Y}_C| < \ell$.
By Claim~\ref{claim:if y' is small, then y^c is big}, we can choose 
$X^* \subset X_C$ and $Y^* \subset Y_C$ such that 
$|X^*| = \lceil n/2 \rceil - \ell - |\overline{Y}_a|$ and 
$|Y^*| = \ell - |\overline{X}_a|$.
In what follows, we find an immersion of $K_{\ell, \lceil n/2 \rceil - \ell}$
with bipartition $(Y^* \cup \overline{X}_a, X^* \cup \overline{Y}_a)$
(see Figure~\ref{fig:case1}).

\begin{figure}[h]
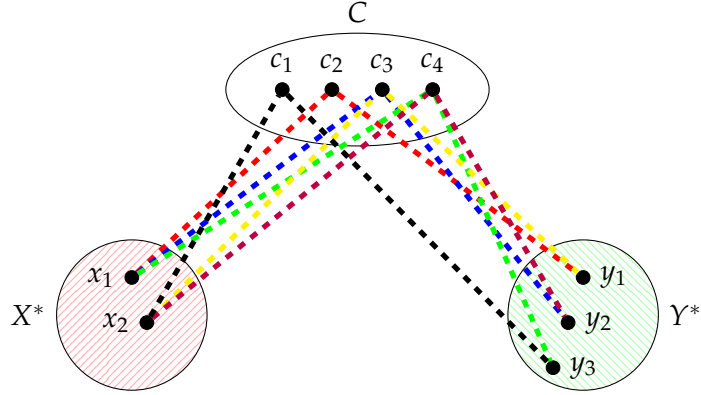

  \centering
  \tikz{
  \node[ellipse,draw,label=above:{$C$},minimum height=1.5cm, minimum width=3.5cm] at (3,4) (C) {};
  \node[ellipse,draw,label=left:{$X^*$},pattern=north east lines,pattern color=red!30,minimum height=2cm, minimum width=2cm] at (0,1) (X) {};
  \node[ellipse,draw,label=right:{$Y^*$},minimum height=2cm, minimum width=2cm,pattern=north west lines,pattern color=green!30] at (6,1) (Y) {};
  \node[circle,draw,fill,label=left:{$x_1$},scale=.5] at (0,1.5) (x1) {};
  \node[circle,draw,fill,label=left:{$x_2$},scale=.5] at (.2,.9) (x2) {};
  \node[circle,draw,fill,label=right:{$y_1$},scale=.5] at (6,1.5) (y1) {};
  \node[circle,draw,fill,label=right:{$y_2$},scale=.5] at (5.8,.9) (y2) {};
  \node[circle,draw,fill,label=right:{$y_3$},scale=.5] at (5.6,.3) (y3) {};

  \node[circle,draw,fill,label=above:{$c_1$},scale=.5] at (2,4) (c1) {};

  \node[circle,draw,fill,label=above:{$c_2$},scale=.5] at (2.66,4) (c2) {};

  \node[circle,draw,fill,label=above:{$c_3$},scale=.5] at (3.33,4) (c3) {};
  \node[circle,draw,fill,label=above:{$c_4$},scale=.5] at (4,4) (c4) {};

  \draw (x1) [dashed,red,line width=2] to (c2) to (y1);
  \draw (x1) [dashed,blue,line width=2] to (c3) to (y2);
  \draw (x1) [dashed,green,line width=2] to (c4) to (y3);
  \draw (x2) [dashed,yellow,line width=2] to (c3) to (y1);
  \draw (x2) [dashed,purple,line width=2] to (c4) to (y2);
  \draw (x2) [dashed,line width=2] to (c1) to (y3);
  }
  \caption{Construction for Case 1, connecting $X^*$ to $Y^*$ through $C$.}
  \label{fig:immersion through c}
\end{figure}

Recall that~$X$, $Y$, and $\overline{X}_a \cup \overline{Y}_a$ each induce 
a clique (each is contained in the non-neighborhood of a vertex and~$G$ has
independence number~2), and thus~$G$ contains all edges joining 
(i) vertices in~$\overline{X}_a$ to vertices in~$\overline{Y}_a$; 
(ii) vertices in~$\overline{X}_a$ to vertices in~$X^*$; and 
(iii) vertices of~$\overline{Y}_a$ to vertices in~$Y^*$.
It remains to find edge-disjoint paths joining vertices in~$X^*$ to vertices
in~$Y^*$.
For that, we use only edges incident to vertices in $C$, which guarantees that our
new paths do not use any previously used edge.
Let $X^* = \{x_1, \ldots, x_{|X^*|}\}$, $Y^* = \{y_1, \ldots, y_{|Y^*|}\}$, and
$C = \{c_1, \ldots, c_{|C|}\}$.
For each pair of vertices $x_i,y_j$, consider the path $x_i c_{i+j} y_j$, where
the additions in the subscript are taken modulo~$|C|$
(see~Figure~\ref{fig:immersion through c}).

It is easy to see that these paths are edge-disjoint if and only if 
$|C| \geq |X^*|,|Y^*|$.
Since, Claim~\ref{claim:bounds_Xa_Ya} gives us $|X_C|, |Y_a| \leq \ell - 2$,
and since in Case 1 we assume $|\overline{X}_C| < \ell$, we have
\begin{align}\label{eq:c-n2-l-ya+3}
  |C| &= n - |\overline{X}_C| - |X_C| - |\overline{Y}_a| - |Y_a|  \nonumber \\ 
      &> n - \ell - (\ell-2) - |\overline{Y}_a| - (\ell-2)   \nonumber\\
      &= n - 3\ell - |\overline{Y}_a| + 4  \nonumber\\
      &= \lceil n/2 \rceil + \lfloor n/2 \rfloor - 3\ell - |\overline{Y}_a| + 4  \nonumber\\
      &\geq \lceil n/2 \rceil + 2\ell - 1 - 3\ell - |\overline{Y}_a| + 4 \qquad  \text{ by~\eqref{eq:size}} \nonumber\\
      &= \lceil n/2 \rceil - \ell - |\overline{Y}_a| + 3  \, .
\end{align} 
Hence, $|C| > \lceil n/2 \rceil - \ell - |\overline{Y}_a| + 3 > |X^*|$.
Analogously to~\eqref{eq:c-n2-l-ya+3}, we can show that
$|C| \geq \lceil n/2 \rceil - \ell - |\overline{X}_a| + 3$, and by~\eqref{eq:size},
that $\lceil n/2 \rceil - \ell - |\overline{X}_a| + 3 > |Y^*|$, as desired.
Since these $X^*,Y^*$-paths are mutually edge-disjoint and edge-disjoint from
all previous paths, we get the desired immersion.

\smallskip\noindent\textbf{Case 2}, 
$|\overline{X}_C| \geq \ell$ or $|\overline{Y}_C| \geq \ell$.
Let \(\gamma = \min\{|\overline{X}_C|,|\overline{Y}_C|\}\).
The proof of this case is again divided into two cases, depending on whether 
\(\gamma \geq \ell\) or \(\gamma < \ell\).
Since the constructions of the immersions in these two cases are similar,
we unify them as follows.
Recall that by Claim~\ref{claim:y-y(a) is big},
either 
  $|X_a| \geq 2 (\ell - |\overline{Y}_a|)$ or 
  $|Y_a| \geq 2 (\ell - |\overline{X}_a|)$. 
Then, if $\gamma\geq \ell$, 
we assume without loss of generality that $|Y_a| \geq 2 (\ell - |\overline{X}_a|)$;
and if $\gamma <\ell$, we assume without loss of generality that $|\overline{Y}_C|=\gamma$ and that $|\overline{X}_C|\geq \ell$.
Notice that in both cases we have \(|\overline{X}_C| \geq \ell\).

\begin{figure}[h]
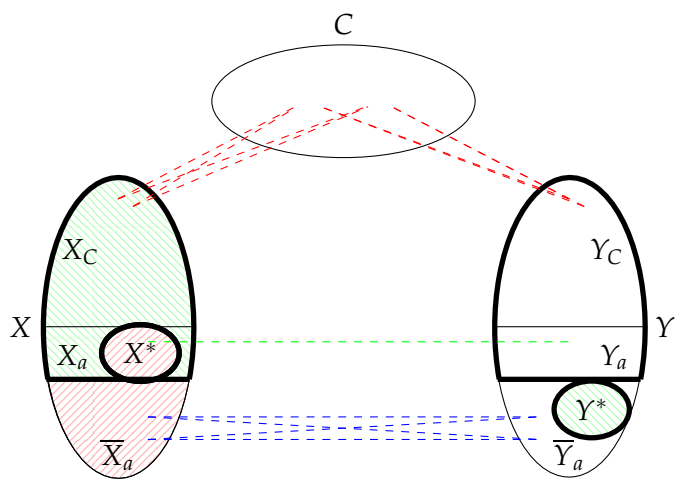

  \centering
  \tikz{
  \node[ellipse,draw,label=above:{$C$},minimum height=1.5cm, minimum width=3.5cm] at (3,3) (C) {};
  \node[ellipse,pattern=north west lines,pattern color=green!30,draw,label=left:{$X$},minimum height=4cm, minimum width=2cm] at (0,0) (X) {};
  \node[ellipse,draw,label=right:{$Y$},minimum height=4cm, minimum width=2cm] at (6,0) (Y) {};
  \draw[fill=white] (0.95,-0.7) arc(-20:-160:1cm and 2cm);
  \draw[pattern=north east lines,pattern color=red!30] (0.95,-0.7) arc(-20:-160:1cm and 2cm);
  \node[ellipse,fill=white,draw,minimum height=0.4cm, minimum width=1cm,line width=2pt] at (0.3,-0.35) (Xs) {\phantom{$X^*$}};
  \node[ellipse,pattern=north east lines,pattern color=red!30,draw,minimum height=0.4cm, minimum width=1cm,line width=2pt] at (0.3,-0.35) (Xs) {$X^*$};
  \node[ellipse,pattern=north west lines,pattern color=green!30,draw,minimum height=0.4cm, minimum width=1cm,line width=2pt] at (6.3,-1.1) (Ys) {$Y^*$};

  \draw[line width=2pt] (0.95,-0.7) arc(-20:200:1cm and 2cm);
  \draw[line width=2pt] (-0.95,-0.7) to (0.95,-0.7);
  \node[circle] at (-0.6, -0.4) (Xa) {$X_a$};

  \draw[line width=2pt] (6.95,-0.7) arc(-20:200:1cm and 2cm);
  \draw[line width=2pt] (5.05,-0.7) to (6.95,-0.7);
  \node[circle] at (6.6, -0.4) (Ya) {$Y_a$};

  \node[circle] at (2.5,3) (a) {};
  \node[circle] at (3.5,3) (b) {};
  \node[circle] at (-.5,1) (XC) {$X_C$};
  \node[circle] at (6.5,1) (YC) {$Y_C$};
  \node[circle] at (0,-1.7) (Xp) {$\overline{X}_a$};
  \node[circle] at (6,-1.7) (Yp) {$\overline{Y}_a$};

  \draw (-1,0) to (1,0);
  \draw (5,0) to (7,0);

  \draw (0,1.7) [dashed, red] to (a);
  \draw (0,1.7) [dashed, red] to (b);
  \draw (0.2,1.6) [dashed, red] to (a);
  \draw (0.2,1.6) [dashed, red] to (b);
  \draw (6,1.7) [dashed, red] to (a);
  \draw (6,1.7) [dashed, red] to (b);
  \draw (6.2,1.6) [dashed, red] to (a);
  \draw (6.2,1.6) [dashed, red] to (b);

  \draw (0.4,-1.5) [dashed, blue] to (5.6,-1.5);
  \draw (0.4,-1.2) [dashed, blue] to (5.6,-1.2);
  \draw (0.4,-1.5) [dashed, blue] to (5.6,-1.2);
  \draw (0.4,-1.2) [dashed, blue] to (5.6,-1.5);

  \draw (0.4,-.2) [dashed, green] to (6,-.2);
  }
  \caption{Immersion for Case 2, when $|\overline{X}_C| \geq \ell$: between $X^* \cup \overline{X}_a$ and $(X_a \setminus X^*) \cup Y^*$.}
  \label{fig:case2}
\end{figure}

By Claim~\ref{claim:bounds_Xa_Ya}, we can choose $Y^* \subset \overline{Y}_a$
with $|Y^*| = \lceil n/2 \rceil - |X|$, and since $|\overline{X}_C| \geq
\ell$, we can choose $X^* \subset \overline{X}_C \setminus \overline{X}_a$
with $|X^*| = \ell - |\overline{X}_a|$.
Note that \(|X_a| - |X^*| + |Y^*| = |X_a| - \ell + |\overline{X}_a| + |Y^*| = |X| + |Y^*|-\ell = \lceil n/2 \rceil -\ell\).
%
Now, we show that
\begin{equation}\label{eq:stars}
    |X^*| \le |Y^*|.
\end{equation} 
Indeed, by Claim~\ref{claim:bounds_Xa_Ya}, we have $|X_a| \leq \ell - 2$, and
hence
\begin{align*}
    |Y^*| &= \lceil n/2 \rceil - |X| \\
          &= \lceil n/2 \rceil - |X_a| - |\overline{X}_a| \\
          &\geq 2\ell - (\ell-2) - |\overline{X}_a| \\
          &= \ell + 2 - |\overline{X}_a| \\
          &> \ell - |\overline{X}_a| \\
          &= |X^*| \, .
\end{align*} 
In what follows, we find an immersion of $K_{\ell, \lceil n/2\rceil - \ell}$ 
with bipartition $\big(X^* \cup \overline{X}_a, (X_a \setminus X^*) \cup Y^*\big)$ 
(see Figure~\ref{fig:case2}).

As before, since $X$, $Y$, and $\overline{X}_a \cup \overline{Y}_a$ induce complete graphs, $G$ contains all edges joining 
(i) vertices in~$\overline{X}_a$ to vertices in~$Y^*$; 
(ii) vertices in~$\overline{X}_a$ to vertices in~$X_a \setminus X^*$; and
(iii) vertices in~$X^*$ to vertices in~$X_a \setminus X^*$.
It remains to find edge-disjoint paths joining vertices in~$X^*$ to vertices
in~$Y^*$.
For these paths, we only use edges that are incident to vertices in~$Y$ and not
to vertices in $\overline{X}_a$; this assures that they are disjoint from the
edges already used.
Let $X^* = \{v_1, v_2, \ldots, v_{|X^*|}\}$ and 
$Y^* = \{y_1, y_2, \ldots, y_{|Y^*|}\}$.
The first step is to use Lemma~\ref{lem:matching_lemma} to find paths joining
each vertex~$v_i$ to all vertices in~$Y^*$ allowing edges between vertices
of~$Y^*$ to be used at most twice.
This is done in the next claim, which is a key step in the proof.

\begin{claim}
\label{claim: construction using matching lemma}
  For each $i \in \{1, 2, \ldots, |X^*|\}$, there is a subgraph $K(v_i)$ which
  contains an immersion of $K_{v_i,Y^*}$ and satisfies that:
  \begin{enumerate}[i)] 
    \item\label{claim:construction:length-at-most-2} each path of $K(v_i)$ with
    an endvertex in~$v_i$ has length at most~$2$;  
    \item\label{claim:construction:z-is-in-Y'} for each path $v_izy_j$ in $K(v_i)$ we have 
        $z \in \overline{Y}_C$; and
    \item\label{claim:construction:multiplicity-at-most-2} if $i \neq j$ and 
        $uw \in E(K(v_i)) \cap E(K(v_j))$, then there is no $r$ with $r \ne i$ and $r \ne j$ such that 
        $uw \in E(K(v_r))$.
        Moreover, if \(v_iz_{i,k}y_{k}\) and \(v_jz_{j,k'}y_{k'}\) 
        are the paths in, respectively, \(K(v_i)\) and \(K(v_j)\) that contain \(uv\), then \(\{u,v\} = \{y_{k},y_{k'}\}\).
  \end{enumerate}
\end{claim}
\begin{proof}
  Note that, since $X^*\subseteq \overline{X}_C$,
  Claim~\ref{claim:neighborhoods_Xprime_Yprime} assures that for each 
  $i \in \{1,2,\ldots, |X^*|\}$, we have 
  $|N(v_i) \cap \overline{Y}_C| > \lceil n/2 \rceil - |X| = |Y^*|$.
  In order to use Lemma~\ref{lem:matching_lemma}, we define, for each such~$i$,
  a set $N_i \subset N(v_i) \cap \overline{Y}_C$ with $|N_i| = |Y^*|$, and also a set of
  auxiliary vertices $A = \{a_1,a_2,\ldots,a_{|Y^*|}\}$ with
  $N(a_j) = \overline{Y}_C$.
  By~\eqref{eq:stars}, we can apply Lemma~\ref{lem:matching_lemma} to
  $N_1,\ldots,N_{|X^*|}$ together with~$A$ to obtain disjoint matchings
  $M_1,M_2,\ldots,M_{|X^*|}$ such that~$M_i$ matches~$A$ to~$N_i$, for
  $i \in \{1, \ldots, |X^*|\}$.
  Let $M_i = \{z_{i,1}a_1, \ldots, z_{i,|Y^*|}a_{|Y^*|}\}$ where $z_{i,j} \in N_i$
  for every $i,j$ (see Figure~\ref{fig:matching construction setup}).

  \begin{figure}[h]
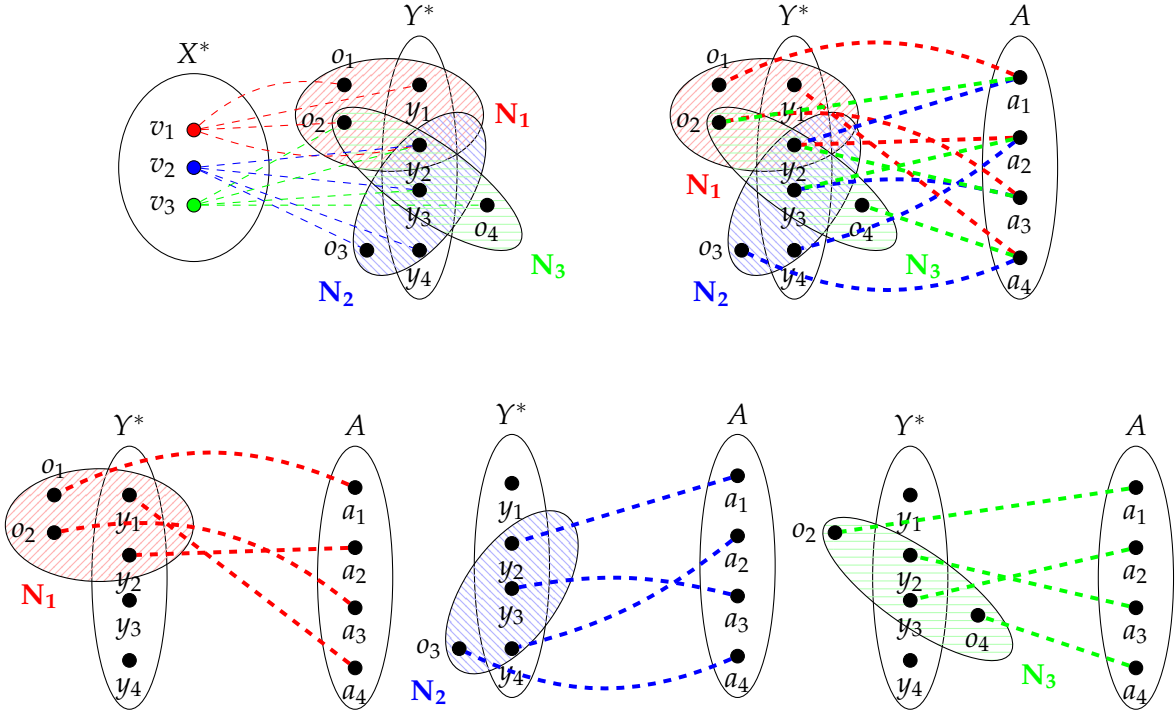

    \centering
    \tikz{
    \begin{pgfonlayer}{background}
      \node[ellipse,draw,label=above:{$X^*$},minimum height=2.5cm, minimum width=2cm] at (-3,0) (X) {};
      \node[ellipse,draw,label=above:{$Y^*$},minimum height=3.5cm, minimum width=1cm] at (0,0) (Y) {};

      \node[ellipse,draw,minimum height=1.5cm, minimum width=2.5cm,pattern=north east lines,pattern color=red!30,label=right:{$\color{red}\mathbf{N_1}$}] at (-.4,.7) (N1) {};
      \node[ellipse,draw,rotate=55,minimum height=1.25cm, minimum width=2.5cm,pattern=north west lines,pattern color=blue!30,label=left:{$\color{blue}\mathbf{N_2}$}] at (0,-.35) (N2) {};
      \node[ellipse,pattern=horizontal lines,rotate=-35,pattern color=green!30,draw,label=right:{$\color{green}\mathbf{N_3}$},minimum height=1cm, minimum width=3cm] at (.1,-.15) (N3) {};
    \end{pgfonlayer}

    \node[circle,draw,fill,scale=.5,label=below:{$y_1$}] at (0,1.1) (y1) {};
    \node[circle,draw,fill,label=below:{$y_2$},scale=.5] at (0,.3) (y2) {};
    \node[circle,draw,fill,label=below:{$y_3$},scale=.5] at (0,-.3) (y3) {};
    \node[circle,draw,fill,label=below:{$y_4$},scale=.5] at (0,-1.1) (y4) {};

    \node[circle,draw,fill=red,scale=.5,label=left:{$v_1$}] at (-3,.5) (x1) {};
    \node[circle,draw,fill=blue,label=left:{$v_2$},scale=.5] at (-3,0) (x2) {};
    \node[circle,draw,fill=green,label=left:{$v_3$},scale=.5] at (-3,-.5) (x3) {};

    \node[circle,draw,fill,scale=.5,label=above:{$o_1$}] at (-1,1.1) (o1) {};
    \node[circle,draw,fill,label=left:{$o_2$},scale=.5] at (-1,.6) (o3) {};
    \node[circle,draw,fill,label=left:{$o_3$},scale=.5] at (-.7,-1.1) (o4) {};
    \node[circle,draw,fill,label=below:{$o_4$},scale=.5] at (.9,-.5) (o5) {};

    \draw (x1) [dashed,red, bend left=25] to (o1);
    \draw (x1) [dashed,red,bend right=15] to (y2);
    \draw (x1) [dashed,red] to (o3);
    \draw (x1) [dashed,red] to (y1);

    \draw (x2) [dashed,blue] to (y2);
    \draw (x2) [dashed,blue] to (y3);
    \draw (x2) [dashed,blue] to (y4);
    \draw (x2) [dashed,blue] to (o4);

    \draw (x3) [dashed,green] to (o3);
    \draw (x3) [dashed,green] to (y2);
    \draw (x3) [dashed,green] to (y3);
    \draw (x3) [dashed,green] to (o5);
    } \hspace{1cm}
    \tikz{
    \begin{pgfonlayer}{background}
      \node[ellipse,draw,label=above:{$Y^*$},minimum height=3.5cm, minimum width=1cm] at (0,0) (Y) {};

      \node[ellipse,draw,minimum height=1.5cm, minimum width=2.5cm,pattern=north east lines,pattern color=red!30,label={[xshift=-0.8cm,yshift=-2cm]$\color{red}\mathbf{N_1}$}] at (-.4,.7) (N1) {};
      \node[ellipse,draw,rotate=55,minimum height=1.25cm, minimum width=2.5cm,pattern=north west lines,pattern color=blue!30,label=left:{$\color{blue}\mathbf{N_2}$}] at (0,-.35) (N2) {};
      \node[ellipse,pattern=horizontal lines,rotate=-35,pattern color=green!30,draw,label=right:{$\color{green}\mathbf{N_3}$},minimum height=1cm, minimum width=3cm] at (.1,-.15) (N3) {};

      \node[ellipse,draw,label=above:{$A$},minimum height=3.5cm, minimum width=1cm] at (3,0) (A) {};
    \end{pgfonlayer}

    \node[circle,draw,fill,scale=.5,label=below:{$y_1$}] at (0,1.1) (y1) {};
    \node[circle,draw,fill,label=below:{$y_2$},scale=.5] at (0,.3) (y2) {};
    \node[circle,draw,fill,label=below:{$y_3$},scale=.5] at (0,-.3) (y3) {};
    \node[circle,draw,fill,label=below:{$y_4$},scale=.5] at (0,-1.1) (y4) {};

    \node[circle,draw,fill,scale=.5,label=below:{$a_1$}] at (3,1.2) (a1) {};
    \node[circle,draw,fill,label=below:{$a_2$},scale=.5] at (3,.4) (a2) {};
    \node[circle,draw,fill,label=below:{$a_3$},scale=.5] at (3,-.4) (a3) {};
    \node[circle,draw,fill,label=below:{$a_4$},scale=.5] at (3,-1.2) (a4) {};

    \node[circle,draw,fill,scale=.5,label=above:{$o_1$}] at (-1,1.1) (o1) {};
    \node[circle,draw,fill,label=left:{$o_2$},scale=.5] at (-1,.6) (o3) {};
    \node[circle,draw,fill,label=left:{$o_3$},scale=.5] at (-.7,-1.1) (o4) {};
    \node[circle,draw,fill,label=below:{$o_4$},scale=.5] at (.9,-.5) (o5) {};

    \draw (a1) [dashed,red,line width=1.5, bend right=25] to (o1);
    \draw (a2) [dashed,red,line width=1.5] to (y2);
    \draw (a4) [dashed,red,line width=1.5] to (y1);
    \draw (a3) [dashed,red,line width=1.5, bend right=25] to (o3);

    \draw (a1) [dashed,blue,line width=1.5] to (y2);
    \draw (a2) [dashed,blue,line width=1.5, bend left=12] to (y4);
    \draw (a3) [dashed,blue,line width=1.5, bend right=12] to (y3);
    \draw (a4) [dashed,blue,line width=1.5,bend left=25] to (o4);
    \draw (a1) [dashed,green,line width=1.5] to (o3);
    \draw (a2) [dashed,green,line width=1.5] to (y3);
    \draw (a3) [dashed,green,line width=1.5] to (y2);
    \draw (a4) [dashed,green,line width=1.5] to (o5);
    }

    \vspace{1cm}
    \tikz{
    \begin{pgfonlayer}{background}
      \node[ellipse,draw,label=above:{$Y^*$},minimum height=3.5cm, minimum width=1cm] at (0,0) (Y) {};

      \node[ellipse,draw,minimum height=1.5cm, minimum width=2.5cm,pattern=north east lines,pattern color=red!30,label={[xshift=-0.8cm,yshift=-2cm]$\color{red}\mathbf{N_1}$}] at (-.4,.7) (N1) {};

      \node[ellipse,draw,label=above:{$A$},minimum height=3.5cm, minimum width=1cm] at (3,0) (A) {};
    \end{pgfonlayer}

    \node[circle,draw,fill,scale=.5,label=below:{$y_1$}] at (0,1.1) (y1) {};
    \node[circle,draw,fill,label=below:{$y_2$},scale=.5] at (0,.3) (y2) {};
    \node[circle,draw,fill,label=below:{$y_3$},scale=.5] at (0,-.3) (y3) {};
    \node[circle,draw,fill,label=below:{$y_4$},scale=.5] at (0,-1.1) (y4) {};

    \node[circle,draw,fill,scale=.5,label=below:{$a_1$}] at (3,1.2) (a1) {};
    \node[circle,draw,fill,label=below:{$a_2$},scale=.5] at (3,.4) (a2) {};
    \node[circle,draw,fill,label=below:{$a_3$},scale=.5] at (3,-.4) (a3) {};
    \node[circle,draw,fill,label=below:{$a_4$},scale=.5] at (3,-1.2) (a4) {};

    \node[circle,draw,fill,scale=.5,label=above:{$o_1$}] at (-1,1.1) (o1) {};
    \node[circle,draw,fill,label=left:{$o_2$},scale=.5] at (-1,.6) (o3) {};

    \draw (a1) [dashed,red,line width=1.5, bend right=25] to (o1);
    \draw (a2) [dashed,red,line width=1.5] to (y2);
    \draw (a4) [dashed,red,line width=1.5] to (y1);
    \draw (a3) [dashed,red,line width=1.5, bend right=25] to (o3);
    }
    \tikz{
    \begin{pgfonlayer}{background}
      \node[ellipse,draw,label=above:{$Y^*$},minimum height=3.5cm, minimum width=1cm] at (0,0) (Y) {};
      \node[ellipse,draw,rotate=55,minimum height=1.25cm, minimum width=2.5cm,pattern=north west lines,pattern color=blue!30,label=left:{$\color{blue}\mathbf{N_2}$}] at (0,-.35) (N2) {};

      \node[ellipse,draw,label=above:{$A$},minimum height=3.5cm, minimum width=1cm] at (3,0) (A) {};
    \end{pgfonlayer}

    \node[circle,draw,fill,scale=.5,label=below:{$y_1$}] at (0,1.1) (y1) {};
    \node[circle,draw,fill,label=below:{$y_2$},scale=.5] at (0,.3) (y2) {};
    \node[circle,draw,fill,label=below:{$y_3$},scale=.5] at (0,-.3) (y3) {};
    \node[circle,draw,fill,label=below:{$y_4$},scale=.5] at (0,-1.1) (y4) {};

    \node[circle,draw,fill,scale=.5,label=below:{$a_1$}] at (3,1.2) (a1) {};
    \node[circle,draw,fill,label=below:{$a_2$},scale=.5] at (3,.4) (a2) {};
    \node[circle,draw,fill,label=below:{$a_3$},scale=.5] at (3,-.4) (a3) {};
    \node[circle,draw,fill,label=below:{$a_4$},scale=.5] at (3,-1.2) (a4) {};

    \node[circle,draw,fill,label=left:{$o_3$},scale=.5] at (-.7,-1.1) (o4) {};

    \draw (a1) [dashed,blue,line width=1.5] to (y2);
    \draw (a2) [dashed,blue,line width=1.5, bend left=12] to (y4);
    \draw (a3) [dashed,blue,line width=1.5, bend right=12] to (y3);
    \draw (a4) [dashed,blue,line width=1.5,bend left=25] to (o4);
    }
    \tikz{
    \begin{pgfonlayer}{background}
      \node[ellipse,draw,label=above:{$Y^*$},minimum height=3.5cm, minimum width=1cm] at (0,0) (Y) {};
      \node[ellipse,pattern=horizontal lines,rotate=-35,pattern color=green!30,draw,label=right:{$\color{green}\mathbf{N_3}$},minimum height=1cm, minimum width=3cm] at (.1,-.15) (N3) {};

      \node[ellipse,draw,label=above:{$A$},minimum height=3.5cm, minimum width=1cm] at (3,0) (A) {};
    \end{pgfonlayer}

    \node[circle,draw,fill,scale=.5,label=below:{$y_1$}] at (0,1.1) (y1) {};
    \node[circle,draw,fill,label=below:{$y_2$},scale=.5] at (0,.3) (y2) {};
    \node[circle,draw,fill,label=below:{$y_3$},scale=.5] at (0,-.3) (y3) {};
    \node[circle,draw,fill,label=below:{$y_4$},scale=.5] at (0,-1.1) (y4) {};

    \node[circle,draw,fill,scale=.5,label=below:{$a_1$}] at (3,1.2) (a1) {};
    \node[circle,draw,fill,label=below:{$a_2$},scale=.5] at (3,.4) (a2) {};
    \node[circle,draw,fill,label=below:{$a_3$},scale=.5] at (3,-.4) (a3) {};
    \node[circle,draw,fill,label=below:{$a_4$},scale=.5] at (3,-1.2) (a4) {};

    \node[circle,draw,fill,label=left:{$o_2$},scale=.5] at (-1,.6) (o3) {};
    \node[circle,draw,fill,label=below:{$o_4$},scale=.5] at (.9,-.5) (o5) {};

    \draw (a1) [dashed,green,line width=1.5] to (o3);
    \draw (a2) [dashed,green,line width=1.5] to (y3);
    \draw (a3) [dashed,green,line width=1.5] to (y2);
    \draw (a4) [dashed,green,line width=1.5] to (o5);
    }
    \caption{The setup for an example of Claim~\ref{claim: construction using matching lemma}.}
    \label{fig:matching construction setup}
  \end{figure}

  For each $v_i \in X^*$, we obtain $K(v_i)$ by using~$y_j$ whenever~$a_j$ is
  used in a matching.
  In other words, for every $1 \leq j \leq |Y^*|$, if $z_{i,j}a_j\in M_i$, then we
  use the path $v_iz_{i,j}y_j$.
  Notice that~$z_{i,j}$ could be~$y_j$ itself.
  When that is the case, we use the path $v_iy_j$.
  Formally, we define 
  \begin{equation*}
    P(i,j) = \begin{cases}
               v_i z_{i,j} y_j & \text{ if } y_j \neq z_{i,j} \\
               v_i y_j & \text{ if } y_j = z_{i,j} \, .
             \end{cases}
  \end{equation*} 
  Notice that $P(i,j)$ may not be edge-disjoint from $P(i,k)$ if $k \neq j$,
  but this can only happen when $P(i,j)=v_iy_ky_j$ and $P(i,k)=v_iy_jy_k$.
  If that is the case, we redefine $P(i,j)$ as $v_iy_j$ and $P(i,k)$ as $v_iy_k$.
  Thus, after doing all the necessary changes, we can assume that $P(i,j)$ is
  disjoint from $P(i,k)$ whenever $j \neq k$ (see Figure~\ref{fig: fixing loops and
    2-cycles}).
  Finally we define $K(v_i) = \bigcup_{j=1}^{|Y^*|} P(i,j)$ and since the
  $P(i,j)$'s are edge disjoint each $K(v_i)$ contains an immersion of
  $K_{v_i,Y^*}$ and clearly satisfies items $i)$ and $ii)$.

  Furthermore, as $M_1,\ldots, M_{|X^*|}$ are disjoint matchings, if 
  $uw \in E(K(v_i)) \cap E(K(v_j))$ for some pair $i \neq j$, then it must be 
  that $u,w \in Y^*$.
  Indeed, by construction \(K(v_i)\) (resp. \(K(v_j)\)) is obtained from a star
  centered in \(v_i\) (resp. \(v_j\)) and with leaves in \(\overline{Y}_C\)
  by adding edges joining its leaves to the vertices in \(Y^*\).
  Then we have \(u,w\notin\{v_i,v_j\}\),
  and hence \(u,w \in \overline{Y}_C\), and \(uw\) 
  is the ``second'' edge of those paths,
  i.e., \(uv\) joins
  a vertex \(z_{i,k_i}\) to \(y_{k_i}\in Y^*\) 
  and a vertex \(z_{j,k_j}\) to \(y_{k_j} \in Y^*\).
  Assume without loss of generality that \(u = z_{i,k_i}\) and \(w = y_{k_i}\).
  Then \(w\in Y^*\).
  Moreover, this implies that \(z_{i,k_i}a_{k_i} \in M_i\).
  If \(u = y_{k_j}\), then \(u\in Y^*\) as desired.
  Otherwise we must have \(u = z_{j,k_j} = z_{i,k_i}\) and \(w = y_{k_j} = y_{k_i}\).
  But this implies that \(z_{i,k_i}a_{k_i} \in M_j\), and hence \(M_i \cap M_j \neq \emptyset\), a contradiction.

  Let $u = y_h$ and $w = y_k$.
  Then either $z_{i,h} = y_k$ or $z_{i,k} = y_h$.
  Assume, without loss of generality, that $z_{i,h} = y_k$.
  This means that $z_{i,h} a_h = y_k a_h \in M_i$.
  Thus $y_k a_h \notin M_r$ for $r \neq i$.
  This, in turn, implies that $z_{j,k} = y_h$, which means that $y_h a_k \in M_j$
  and $y_h a_k \notin M_r$ for $r \neq j$.
  This proves $iii)$.
  \begin{figure}
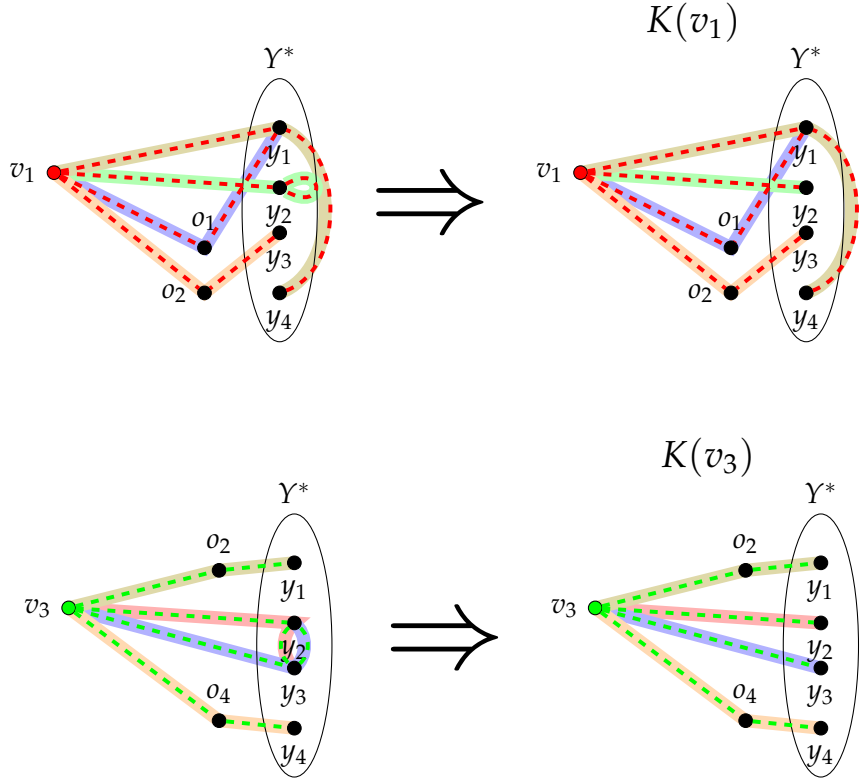

    \centering
    \tikz{
    \node[ellipse,draw,label=above:{$Y^*$},minimum height=3.5cm, minimum width=1cm] at (0,0) (Y) {};
    \node[circle,draw,fill,scale=.5,label=below:{$y_1$}] at (0,1.1) (y1) {};

    \node[circle,draw,fill,label=below:{$y_2$},scale=.5] at (0,.3) (y2) {};
    \node[circle,draw,fill,label=below:{$y_3$},scale=.5] at (0,-.3) (y3) {};
    \node[circle,draw,fill,label=below:{$y_4$},scale=.5] at (0,-1.1) (y4) {};

    \node[circle,draw,fill,label=above:{$o_1$},scale=.5] at (-1,-.5) (o2) {};
    \node[circle,draw,fill,label=left:{$o_2$},scale=.5] at (-1,-1.1) (o3) {};

    \node[circle,draw,fill=red,scale=.5,label=left:{$v_1$}] at (-3,.5) (x1) {};

    \draw (x1) [dashed,red,line width=1.5] to (y2);
    \draw (x1) [dashed,red,line width=1.5] to (o2);
    \draw (x1) [dashed,red,line width=1.5] to (o3);
    \draw (x1) [dashed,red,line width=1.5] to (y1);

    \draw (y1) [dashed,red,line width=1.5] to (o2);
    \draw (y2) [dashed,red,line width=1.5] to  [out=30,in=330,looseness=17] (y2);
    \draw (y3) [dashed,red,line width=1.5] to (o3);
    \draw (y4) [dashed,red,line width=1.5, bend right=70] to (y1);

    \begin{pgfonlayer}{background}
      \draw [line width=5, color=blue!30](x1.center) to (o2.center) to (y1.center);
      \draw [line width=5, color=green!30](x1.center) to (y2.center);
      \draw [line width=5, color=olive!30](x1.center) to (y1.center) [bend left=70] to (y4.center);
      \draw [line width=5, color=orange!30](x1.center) to (o3.center) to  (y3.center);
      \draw (y2) [line width=5, color=green!30] to [out=30,in=330,looseness=17] (y2);
    \end{pgfonlayer}

    \node[circle,scale=4] at (2,0) (m1) {$\Rightarrow$};


    \node[ellipse,draw,label=above:{$Y^*$},minimum height=3.5cm, minimum width=1cm] at (7,0) (Y) {};
    \node[circle,draw,fill,scale=.5,label=below:{$y_1$}] at (7,1.1) (y1) {};

    \node[circle,draw,fill,label=below:{$y_2$},scale=.5] at (7,.3) (y2) {};
    \node[circle,draw,fill,label=below:{$y_3$},scale=.5] at (7,-.3) (y3) {};
    \node[circle,draw,fill,label=below:{$y_4$},scale=.5] at (7,-1.1) (y4) {};

    \node[circle,draw,fill,label=above:{$o_1$},scale=.5] at (6,-.5) (o2) {};
    \node[circle,draw,fill,label=left:{$o_2$},scale=.5] at (6,-1.1) (o3) {};

    \node[circle,draw,fill=red,scale=.5,label=left:{$v_1$}] at (4,.5) (x1) {};

    \draw (x1) [dashed,red,line width=1.5] to (y2);
    \draw (x1) [dashed,red,line width=1.5] to (o2);
    \draw (x1) [dashed,red,line width=1.5] to (o3);
    \draw (x1) [dashed,red,line width=1.5] to (y1);

    \draw (y1) [dashed,red,line width=1.5] to (o2);
    \draw (y3) [dashed,red,line width=1.5] to (o3);
    \draw (y4) [dashed,red,line width=1.5, bend right=70] to (y1);

    \begin{pgfonlayer}{background}
      \draw [line width=5, color=blue!30](x1.center) to (o2.center) to (y1.center);
      \draw [line width=5, color=green!30](x1.center) to (y2.center);
      \draw [line width=5, color=olive!30](x1.center) to (y1.center) [bend left=70] to (y4.center);
      \draw [line width=5, color=orange!30](x1.center) to (o3.center) to  (y3.center);
    \end{pgfonlayer}

    \node[circle,scale=1.3] at (5.5,2.5) (m1) {$K(v_1)$};

    }\\[0.6cm]
    \tikz{

    \node[ellipse,draw,label=above:{$Y^*$},minimum height=3.5cm, minimum width=1cm] at (0,0) (Y) {};
    \node[circle,draw,fill,scale=.5,label=below:{$y_1$}] at (0,1.1) (y1) {};

    \node[circle,draw,fill,label=below:{$y_2$},scale=.5] at (0,.3) (y2) {};
    \node[circle,draw,fill,label=below:{$y_3$},scale=.5] at (0,-.3) (y3) {};
    \node[circle,draw,fill,label=below:{$y_4$},scale=.5] at (0,-1.1) (y4) {};

    \node[circle,draw,fill,label=above:{$o_2$},scale=.5] at (-1,1) (o3) {};
    \node[circle,draw,fill,label=above:{$o_4$},scale=.5] at (-1,-1) (o5) {};

    \node[circle,draw,fill=green,scale=.5,label=left:{$v_3$}] at (-3,.5) (x2) {};

    \draw (x2) [dashed,green,line width=1.5] to (y2);
    \draw (x2) [dashed,green,line width=1.5] to (y3);
    \draw (x2) [dashed,green,line width=1.5] to (o3);
    \draw (x2) [dashed,green,line width=1.5] to (o5);

    \draw (o5) [dashed,green,line width=1.5] to (y4);
    \draw (y2) [dashed,green,line width=1.5,bend right=50] to (y3);
    \draw (y3) [dashed,green,line width=1.5, bend right=50] to (y2);
    \draw (o3) [dashed,green,line width=1.5] to (y1);

    \begin{pgfonlayer}{background}
      \draw [line width=5, color=red!30](x2.center) to (y2.center) to [bend right=50](y3.center);
      \draw [line width=5, color=blue!30](x2.center) to (y3.center) to [bend right=50](y2.center);
      \draw [line width=5, color=olive!30](x2.center) to (o3.center) to (y1.center);
      \draw [line width=5, color=orange!30](x2.center) to (o5.center) to  (y4.center);
    \end{pgfonlayer}

    \node[circle,scale=4] at (2,0) (m1) {$\Rightarrow$};


    \node[ellipse,draw,label=above:{$Y^*$},minimum height=3.5cm, minimum width=1cm] at (7,0) (Y) {};
    \node[circle,draw,fill,scale=.5,label=below:{$y_1$}] at (7,1.1) (y1) {};

    \node[circle,draw,fill,label=below:{$y_2$},scale=.5] at (7,.3) (y2) {};
    \node[circle,draw,fill,label=below:{$y_3$},scale=.5] at (7,-.3) (y3) {};
    \node[circle,draw,fill,label=below:{$y_4$},scale=.5] at (7,-1.1) (y4) {};

    \node[circle,draw,fill,label=above:{$o_2$},scale=.5] at (6,1) (o3) {};
    \node[circle,draw,fill,label=above:{$o_4$},scale=.5] at (6,-1) (o5) {};

    \node[circle,draw,fill=green,scale=.5,label=left:{$v_3$}] at (4,.5) (x2) {};

    \draw (x2) [dashed,green,line width=1.5] to (y2);
    \draw (x2) [dashed,green,line width=1.5] to (y3);
    \draw (x2) [dashed,green,line width=1.5] to (o3);
    \draw (x2) [dashed,green,line width=1.5] to (o5);

    \draw (o5) [dashed,green,line width=1.5] to (y4);
    \draw (o3) [dashed,green,line width=1.5] to (y1);

    \begin{pgfonlayer}{background}
      \draw [line width=5, color=red!30](x2.center) to (y2.center);
      \draw [line width=5, color=blue!30](x2.center) to (y3.center);
      \draw [line width=5, color=olive!30](x2.center) to (o3.center) to (y1.center);
      \draw [line width=5, color=orange!30](x2.center) to (o5.center) to  (y4.center);
    \end{pgfonlayer}

    \node[circle,scale=1.3] at (5.5,2.5) (m1) {$K(v_3)$};
    }
    \caption{Defining $K(v_1)$ and $K(v_3)$ for the graph in Figure~\ref{fig:matching construction setup}.}
    \label{fig: fixing loops and 2-cycles}
  \end{figure}
\end{proof}

From now on, let $K(v_1), \ldots, K(v_{|X^*|})$ be the subgraphs given by
Claim~\ref{claim: construction using matching lemma}.
We would like the $v_i,Y^*$-paths on these subgraphs to be the $X^*,Y^*$-paths
in our immersion.
Yet, if $i \ne j$, $K(v_i)$ might not be edge disjoint from $K(v_j)$.
Fortunately, by Claim~\ref{claim: construction using matching lemma} $iii)$ these
intersections are restricted, and, in what follows, we fix them.

Let~$H$ be the (multi)graph with vertex set~$V(G)$, 
and whose edge set is the disjoint union of $E(K(v_1)), \ldots, E(K(v_{|X^*|}))$, 
where the multiplicity of an edge \(e\) is \(|\{i : e \in E(K(v_i))\}|\).
Notice that by Claim~\ref{claim: construction using matching lemma} $iii)$, the
multiplicity of every edge of~$H$ is at most~2, and the only edges with multiplicity~2 are
between vertices of~$Y^*$.
To fix these repetitions, we take detours through other vertices in~$Y$.
We now need to divide Case 2 into two subcases, depending on whether
$|\overline{Y}_C| \geq \ell$.

\smallskip
\noindent
\textbf{Case 2.1}, $\gamma \geq \ell$.
Recall that in this case, we have $|\overline{X}_C|,|\overline{Y}_C| \geq \ell$, and that we assume, without loss of generality, that
\begin{equation}\label{eq:size_barYa}
    |Y\setminus\overline{Y}_a| = |Y_a| \geq 2(\ell-|\overline{X}_a|) \, .
\end{equation}

Let \(u\in Y^*\).
By Claim~\ref{claim: construction using matching lemma} $iii)$ we know that if~$u$
is incident to an edge of multiplicity~2, say $uw$, then there is a
path in some $K(v_i)$ that contains~$uw$ and ends at~$u$.
Recall that $m_H(e)$ is the multiplicity of~$e$ in~$H$,
and, for a given vertex $u \in Y^*$, let~$p_u$ be defined as 
\[p_u = \sum_{uv\in E(H),\newline v \in Y^*} (m_H(uv) - 1)\, ,\]
and let $q_u = \left|N_{H}(u) \cap \left( Y \setminus Y^* \right) \right|$.
Let us see that $q_u + p_u \leq |X^*| = \ell - |\overline{X}_a|$.
Indeed, there are precisely~$|X^*|$ paths in \(H\) that have~$u$ as an endpoint, where~$p_u$ of these
paths use edges with multiplicity~2,
and~$q_u$ use edges that join \(u\) to vertices in~$Y\setminus Y^*$.
Since edges with multiplicity~2 are incident only to vertices in \(Y^*\),
none of these paths were counted twice,
but there could be some of the \(|X^*|\) paths that are not counted above.
Now we show that if~$uw$ is an edge with multiplicity~2 in~$H$, then
\begin{equation}
\label{eq:qq}
  |Y \setminus Y^*| > q_u + q_w \, .   
\end{equation} 
As $uw$ has multiplicity~2, we have $p_u, p_w \geq 1$.
So using~\eqref{eq:size_barYa} and that $Y^*\subseteq \overline{Y}_a$, we
indeed have 
\begin{align*}
    q_u + q_w &\leq \ell - |\overline{X}_a| - p_u + \ell - |\overline{X}_a| - p_w \\
              &\leq \ell - |\overline{X}_a| - 1 + \ell - |\overline{X}_a| - 1 \\
              &= 2(\ell - |\overline{X}_a| - 1) \\
              &\leq |Y \setminus \overline{Y}_a| - 2  \\
              &\leq |Y \setminus Y^*| - 1 \, .
\end{align*}

We now describe a process that at each step picks an edge~$uw$ with multiplicity \(2\) of~$H$,
reduces~$p_u$ and~$p_w$ by~$1$, increases~$q_u$ and~$q_w$ by~$1$, and
leaves~$p_t$ and~$q_t$ unchanged for every vertex $t \in Y^* \setminus
\{u,w\}$.
In other words, the function $I \colon Y^* \to \mathbb{N}$ such that $I(u) = p_u + q_u$
is an invariant of this process, while the number of edges with multiplicity~2 decreases.

\begin{figure}
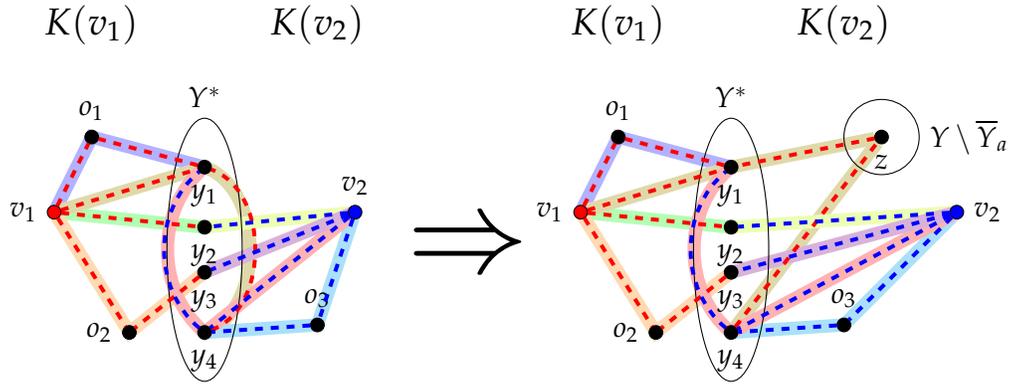

  \centering
  \tikz{
  \node[ellipse,draw,label=above:{$Y^*$},minimum height=3.5cm, minimum width=1cm] at (0,0) (Y) {};
  \node[circle,draw,fill,scale=.5,label=below:{$y_1$}] at (0,1.1) (y1) {};

  \node[circle,draw,fill,label=below:{$y_2$},scale=.5] at (0,.3) (y2) {};
  \node[circle,draw,fill,label=below:{$y_3$},scale=.5] at (0,-.3) (y3) {};
  \node[circle,draw,fill,label=below:{$y_4$},scale=.5] at (0,-1.1) (y4) {};

  \node[circle,draw,fill,label=above:{$o_1$},scale=.5] at (-1.5,1.5) (o2) {};
  \node[circle,draw,fill,label=left:{$o_2$},scale=.5] at (-1,-1.1) (o3) {};

  \node[circle,draw,fill=red,scale=.5,label=left:{$v_1$}] at (-2,.5) (x1) {};

  \draw (x1) [dashed,red,line width=1.5] to (y2);
  \draw (x1) [dashed,red,line width=1.5] to (o2);
  \draw (x1) [dashed,red,line width=1.5] to (o3);
  \draw (x1) [dashed,red,line width=1.5] to (y1);

  \draw (y1) [dashed,red,line width=1.5] to (o2);
  \draw (y3) [dashed,red,line width=1.5] to (o3);
  \draw (y4) [dashed,red,line width=1.5, bend right=70] to (y1);

  \begin{pgfonlayer}{background}
    \draw [line width=5, color=blue!30](x1.center) to (o2.center) to (y1.center);
    \draw [line width=5, color=green!30](x1.center) to (y2.center);
    \draw [line width=5, color=olive!30](x1.center) to (y1.center) [bend left=70] to (y4.center);
    \draw [line width=5, color=orange!30](x1.center) to (o3.center) to  (y3.center);
  \end{pgfonlayer}

  \node[circle,draw,fill,label=above:{$o_3$},scale=.5] at (1.5,-1) (o4) {};
  \node[circle,draw,fill=blue,scale=.5,label=above:{$v_2$}] at (2,.5)(x2) {};
  \draw (x2) [dashed,blue,line width=1.5] to (o4);
  \draw (x2) [dashed,blue,line width=1.5] to (y2);
  \draw (x2) [dashed,blue,line width=1.5] to (y3);
  \draw (x2) [dashed,blue,line width=1.5] to (y4);

  \draw (y4) [dashed,blue,line width=1.5, bend left=50] to (y1);
  \draw (o4) [dashed,blue,line width=1.5] to (y4);

  \begin{pgfonlayer}{background}
    \draw [line width=5, color=lime!30](x2.center) to (y2.center);
    \draw [line width=5, color=cyan!30](x2.center) to (o4.center) to (y4.center);
    \draw [line width=5, color=red!30](x2.center) to (y4.center) [bend left=50] to  (y1.center);
    \draw [line width=5, color=violet!30](x2.center) to (y3.center);
  \end{pgfonlayer}
  \node[circle,scale=1.3] at (-1.5,3) (m1) {$K(v_1)$};
  \node[circle,scale=1.3] at (1.5,3) (m1) {$K(v_2)$};
  \node[circle,scale=4] at (3.5,0) (m1) {$\Rightarrow$};


  \node[ellipse,draw,label=above:{$Y^*$},minimum height=3.5cm, minimum width=1cm] at (7,0) (Y) {};

  \node[ellipse,draw,label=right:{$Y\setminus \overline{Y}_a$},minimum height=1cm, minimum width=1cm] at (9,1.5) (Y) {};
  \node[circle,draw,fill,scale=.5,label=below:{$y_1$}] at (7,1.1) (y1) {};
  \node[circle,draw,fill,scale=.5,label=below:{$z$}] at (9,1.5) (z) {};
  \node[circle,draw,fill,label=below:{$y_2$},scale=.5] at (7,.3) (y2) {};
  \node[circle,draw,fill,label=below:{$y_3$},scale=.5] at (7,-.3) (y3) {};
  \node[circle,draw,fill,label=below:{$y_4$},scale=.5] at (7,-1.1) (y4) {};

  \node[circle,draw,fill,label=above:{$o_1$},scale=.5] at (5.5,1.5) (o2) {};
  \node[circle,draw,fill,label=left:{$o_2$},scale=.5] at (6,-1.1) (o3) {};

  \node[circle,draw,fill=red,scale=.5,label=left:{$v_1$}] at (5,.5) (x1) {};

  \draw (x1) [dashed,red,line width=1.5] to (y2);
  \draw (x1) [dashed,red,line width=1.5] to (o2);
  \draw (x1) [dashed,red,line width=1.5] to (o3);
  \draw (x1) [dashed,red,line width=1.5] to (y1);

  \draw (y1) [dashed,red,line width=1.5] to (o2);
  \draw (y3) [dashed,red,line width=1.5] to (o3);
  \draw (y4) [dashed,red,line width=1.5] to (z) to (y1);

  \begin{pgfonlayer}{background}
    \draw [line width=5, color=blue!30](x1.center) to (o2.center) to (y1.center);
    \draw [line width=5, color=green!30](x1.center) to (y2.center);
    \draw [line width=5, color=olive!30](x1.center) to (y1.center) to (z) to (y4.center);
    \draw [line width=5, color=orange!30](x1.center) to (o3.center) to  (y3.center);
  \end{pgfonlayer}

  \node[circle,draw,fill,label=above:{$o_3$},scale=.5] at (8.5,-1) (o4) {};
  \node[circle,draw,fill=blue,scale=.5,label=right:{$v_2$}] at (10,.5)(x2) {};
  \draw (x2) [dashed,blue,line width=1.5] to (o4);
  \draw (x2) [dashed,blue,line width=1.5] to (y2);
  \draw (x2) [dashed,blue,line width=1.5] to (y3);
  \draw (x2) [dashed,blue,line width=1.5] to (y4);

  \draw (y4) [dashed,blue,line width=1.5, bend left=50] to (y1);
  \draw (o4) [dashed,blue,line width=1.5] to (y4);

  \begin{pgfonlayer}{background}
    \draw [line width=5, color=lime!30](x2.center) to (y2.center);
    \draw [line width=5, color=cyan!30](x2.center) to (o4.center) to (y4.center);
    \draw [line width=5, color=red!30](x2.center) to (y4.center) [bend left=50] to  (y1.center);
    \draw [line width=5, color=violet!30](x2.center) to (y3.center);
  \end{pgfonlayer}
  \node[circle,scale=1.3] at (5.5,3) (m1) {$K(v_1)$};
  \node[circle,scale=1.3] at (8.5,3) (m1) {$K(v_2)$};
  }
  \caption{Removing a multiple edge through $Y\setminus \overline{Y}_a$.}
  \label{fig: removing a double edge}
\end{figure}

Let $v_1,v_2 \in X^*$ be such that $uw \in E(K(v_1)) \cap E(K(v_2))$.
By~\eqref{eq:qq}, there is a vertex $z \in Y \setminus Y^*$ such that~$z$ is 
neither adjacent to~$u$ nor~$w$ in~$H$.
Modify~$H$ by replacing the path $v_1 u w$ in $K(v_1)$ with the path $v_1 u z w$,
and leave $v_2wu$ unchanged (see Figure~\ref{fig: removing a double edge}).
Then this replacement does not increase the multiplicity of any edge, reduces
the multiplicity of edge $uv$ by~$1$, and increases the sizes of 
$N_H(u) \cap Y \setminus Y^*$ and $N_H(w) \cap Y \setminus Y^*$ by~$1$,
because after switching we have \(z\in N_H(u)\cap N_H(w)\).
Furthermore, as $z \in Y \setminus Y^*$, it does not change the values of~$p_t$
and~$q_t$ for any vertex $t \in Y^* \setminus \{u,w\}$.
Hence after each step we still have 
$q_u + p_u \leq |X^*| = \ell - |\overline{X}_a|$ for every $u \in Y^*$, and
thus we can continue to assume~\eqref{eq:qq}.
Repeating this process until $H$ becomes a simple graph ends up giving us that~$H$
contains an immersion of~$K_{X^*,Y^*}$ and therefore~$G$ contains an
immersion of $K_{\ell,\lceil n/2 \rceil -\ell}$.

\smallskip
\noindent
\textbf{Case 2.2, $\gamma < \ell$.}
Recall that in this case we assume $\gamma = |\overline{Y}_C| < \ell$.
In this case, we modify the subgraphs given by Claim~\ref{claim: construction using
  matching lemma} so that the set of edges of each $K(v_i)$ joining two vertices in \(Y^*\) induces a matching.
For each $1 \leq i \leq |X^*|$ let $Q(v_i)$ be the subgraph of $K(v_i)$ induced
by the vertices of~$Y^*$.
Notice that for each vertex $y_j \in Y^*$, the degree of $y_j$ in $Q(v_i)$ is
the number of times~$y_j$ is in the middle of some path $v_iz_{i,k}y_k$
(i.e., $y_j = z_{i,k}$), which is at most~1 because~$v_iy_j$ can be
used only once, plus the number of times~$y_j$ is the end of some path
$v_iz_{i,k}y_k$, which is~1 because there is precisely one path in
$K(v_i)$ joining~$v_i$ to~$y_j$.
Thus we have $d_{Q(v_i)}(y_j) \leq 2$, and so each component of~$Q(v_i)$ is
either a path or a cycle.

\begin{figure}[h]
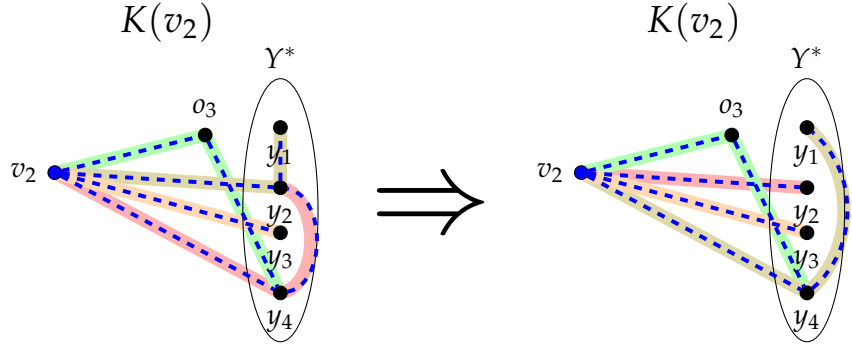

  \centering
  \tikz{

  \node[ellipse,draw,label=above:{$Y^*$},minimum height=3.5cm, minimum width=1cm] at (0,0) (Y) {};
  \node[circle,draw,fill,scale=.5,label=below:{$y_1$}] at (0,1.1) (y1) {};
  \node[circle,draw,fill,label=below:{$y_2$},scale=.5] at (0,.3) (y2) {};
  \node[circle,draw,fill,label=below:{$y_3$},scale=.5] at (0,-.3) (y3) {};
  \node[circle,draw,fill,label=below:{$y_4$},scale=.5] at (0,-1.1) (y4) {};
  \node[circle,draw,fill,label=above:{$o_3$},scale=.5] at (-1,1) (o4) {};
  \node[circle,draw,fill=blue,scale=.5,label=left:{$v_2$}] at (-3,.5) (x2) {};

  \draw (x2) [dashed,blue,line width=1.5] to (o4);
  \draw (x2) [dashed,blue,line width=1.5] to (y2);
  \draw (x2) [dashed,blue,line width=1.5] to (y3);
  \draw (x2) [dashed,blue,line width=1.5] to (y4);

  \draw (y4) [dashed,blue,line width=1.5, bend right=80] to (y2);
  \draw (o4) [dashed,blue,line width=1.5] to (y4);
  \draw (y2) [dashed,blue,line width=1.5] to (y1);

  \begin{pgfonlayer}{background}
    \draw [line width=5, color=red!30](x2.center) to (y4.center) [bend right=80] to (y2.center);
    \draw [line width=5, color=green!30](x2.center) to (o4.center) to (y4.center);
    \draw [line width=5, color=olive!30](x2.center) to (y2.center) to  (y1.center);
    \draw [line width=5, color=orange!30](x2.center) to (y3.center);
  \end{pgfonlayer}
  \node[circle,scale=1.3] at (-1.5,2.5) (m1) {$K(v_2)$};

  \node[circle,scale=4] at (2,0) (m1) {$\Rightarrow$};

  \node[ellipse,draw,label=above:{$Y^*$},minimum height=3.5cm, minimum width=1cm] at (7,0) (Y) {};
  \node[circle,draw,fill,scale=.5,label=below:{$y_1$}] at (7,1.1) (y1) {};
  \node[circle,draw,fill,label=below:{$y_2$},scale=.5] at (7,.3) (y2) {};
  \node[circle,draw,fill,label=below:{$y_3$},scale=.5] at (7,-.3) (y3) {};
  \node[circle,draw,fill,label=below:{$y_4$},scale=.5] at (7,-1.1) (y4) {};
  \node[circle,draw,fill,label=above:{$o_3$},scale=.5] at (6,1) (o4) {};
  \node[circle,draw,fill=blue,scale=.5,label=left:{$v_2$}] at (4,.5) (x2) {};

  \draw (x2) [dashed,blue,line width=1.5] to (o4);
  \draw (x2) [dashed,blue,line width=1.5] to (y2);
  \draw (x2) [dashed,blue,line width=1.5] to (y3);
  \draw (x2) [dashed,blue,line width=1.5] to (y4);

  \draw (y4) [dashed,blue,line width=1.5, bend right=50] to (y1);
  \draw (o4) [dashed,blue,line width=1.5] to (y4);

  \begin{pgfonlayer}{background}
    \draw [line width=5, color=red!30](x2.center) to (y2.center);
    \draw [line width=5, color=green!30](x2.center) to (o4.center) to (y4.center);
    \draw [line width=5, color=olive!30](x2.center) to (y4.center) [bend right=50] to  (y1.center);
    \draw [line width=5, color=orange!30](x2.center) to (y3.center);
  \end{pgfonlayer}
  \node[circle,scale=1.3] at (5.5,2.5) (m1) {$K(v_2)$};
  }
  \caption{Removing a path in $Y^*$.}
  \label{fig: removing paths}
\end{figure}

If $t_1 t_2 \cdots t_s t_1$ is a cycle in $Q(v_i)$, then we replace, in
$K(v_i)$, the path $v_i t_{j} t_{j+1}$ with the path $v_i t_{j+1}$, for every
$1 \leq j \leq s$, where $t_{s+1} = t_1$.
If $t_1 t_2 \cdots t_s$ is a component of $Q(v_i)$ that is a path, then we
replace, in $K(v_i)$, the path $v_i t_{j} t_{j+1}$ with the path $v_i t_{j+1}$,
for each $1 \leq j \leq s-2$, and replace the path $v_i t_{s-1} t_s$ with the
path $v_i t_1 t_s$ (see Figure~\ref{fig: removing paths}).
After these operations $E(Q(v_i))$ is a matching.

\begin{figure}[h]
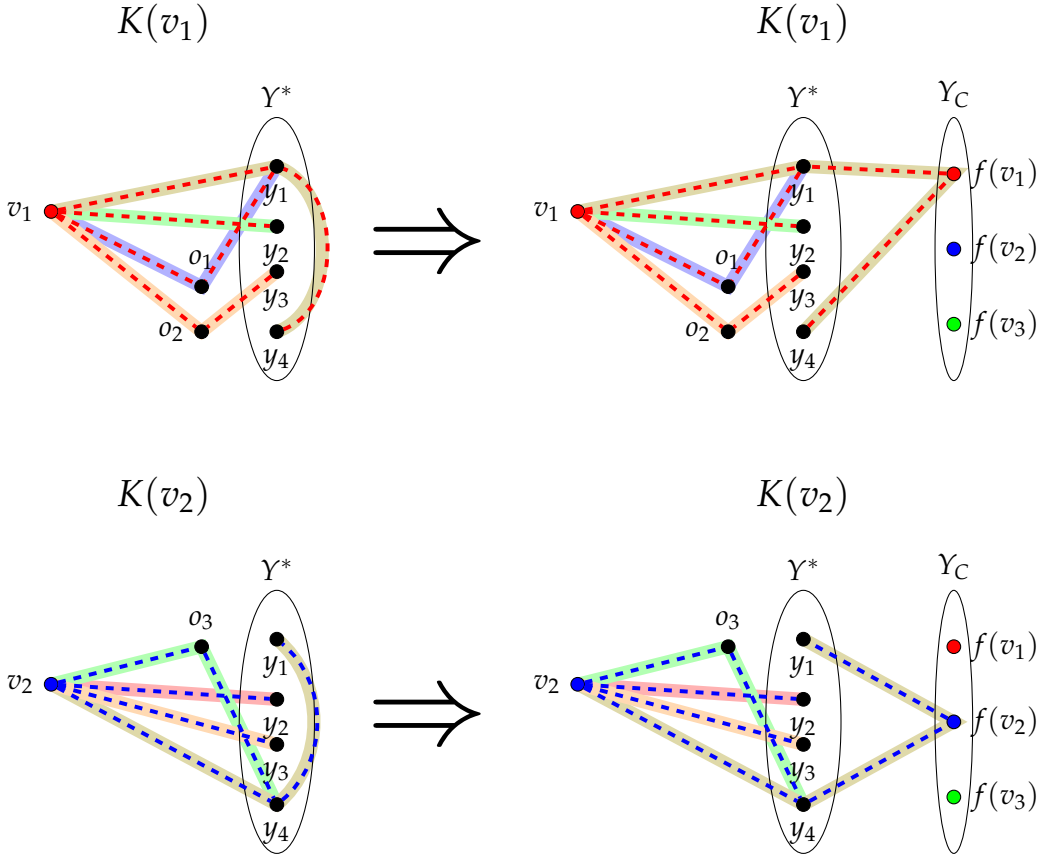

  \centering
  \tikz{
  \node[ellipse,draw,label=above:{$Y^*$},minimum height=3.5cm, minimum width=1cm] at (0,0) (Y) {};
  \node[circle,draw,fill,scale=.5,label=below:{$y_1$}] at (0,1.1) (y1) {};

  \node[circle,draw,fill,label=below:{$y_2$},scale=.5] at (0,.3) (y2) {};
  \node[circle,draw,fill,label=below:{$y_3$},scale=.5] at (0,-.3) (y3) {};
  \node[circle,draw,fill,label=below:{$y_4$},scale=.5] at (0,-1.1) (y4) {};

  \node[circle,draw,fill,label=above:{$o_1$},scale=.5] at (-1,-.5) (o2) {};
  \node[circle,draw,fill,label=left:{$o_2$},scale=.5] at (-1,-1.1) (o3) {};

  \node[circle,draw,fill=red,scale=.5,label=left:{$v_1$}] at (-3,.5) (x1) {};

  \draw (x1) [dashed,red,line width=1.5] to (y2);
  \draw (x1) [dashed,red,line width=1.5] to (o2);
  \draw (x1) [dashed,red,line width=1.5] to (o3);
  \draw (x1) [dashed,red,line width=1.5] to (y1);

  \draw (y1) [dashed,red,line width=1.5] to (o2);
  \draw (y3) [dashed,red,line width=1.5] to (o3);
  \draw (y4) [dashed,red,line width=1.5, bend right=70] to (y1);

  \begin{pgfonlayer}{background}
    \draw [line width=5, color=blue!30](x1.center) to (o2.center) to (y1.center);
    \draw [line width=5, color=green!30](x1.center) to (y2.center);
    \draw [line width=5, color=olive!30](x1.center) to (y1.center) [bend left=70] to (y4.center);
    \draw [line width=5, color=orange!30](x1.center) to (o3.center) to  (y3.center);
  \end{pgfonlayer}

  \node[circle,scale=1.3] at (-1.5,3) (m1) {$K(v_1)$};
  \node[circle,scale=4] at (2,0) (m1) {$\Rightarrow$};


  \node[ellipse,draw,label=above:{$Y^*$},minimum height=3.5cm, minimum width=1cm] at (7,0) (Y) {};

  \node[ellipse,draw,label=above:{$Y_C$},minimum height=3.5cm, minimum width=.5cm] at (9,0) (Y) {};
  \node[circle,draw,fill=red,scale=.5,label=right:{$f(v_1)$}] at (9,1) (fv1) {};
  \node[circle,draw,fill=blue,scale=.5,label=right:{$f(v_2)$}] at (9,0) (fv2) {};
  \node[circle,draw,fill=green,scale=.5,label=right:{$f(v_3)$}] at (9,-1) (fv3) {};

  \node[circle,draw,fill,scale=.5,label=below:{$y_1$}] at (7,1.1) (y1) {};

  \node[circle,draw,fill,label=below:{$y_2$},scale=.5] at (7,.3) (y2) {};
  \node[circle,draw,fill,label=below:{$y_3$},scale=.5] at (7,-.3) (y3) {};
  \node[circle,draw,fill,label=below:{$y_4$},scale=.5] at (7,-1.1) (y4) {};

  \node[circle,draw,fill,label=above:{$o_1$},scale=.5] at (6,-.5) (o2) {};
  \node[circle,draw,fill,label=left:{$o_2$},scale=.5] at (6,-1.1) (o3) {};

  \node[circle,draw,fill=red,scale=.5,label=left:{$v_1$}] at (4,.5) (x1) {};

  \draw (x1) [dashed,red,line width=1.5] to (y2);
  \draw (x1) [dashed,red,line width=1.5] to (o2);
  \draw (x1) [dashed,red,line width=1.5] to (o3);
  \draw (x1) [dashed,red,line width=1.5] to (y1);

  \draw (y1) [dashed,red,line width=1.5] to (o2);
  \draw (y3) [dashed,red,line width=1.5] to (o3);
  \draw (y4) [dashed,red,line width=1.5]
  to (fv1) to (y1);

  \begin{pgfonlayer}{background}
    \draw [line width=5, color=blue!30](x1.center) to (o2.center) to (y1.center);
    \draw [line width=5, color=green!30](x1.center) to (y2.center);
    \draw [line width=5, color=olive!30](x1.center) to (y1.center) to (fv1.center) to (y4.center);
    \draw [line width=5, color=orange!30](x1.center) to (o3.center) to  (y3.center);
  \end{pgfonlayer}

  \node[circle,scale=1.3] at (7,3) (m1) {$K(v_1)$};

  }\\[0.6cm]
  \tikz{

  \node[ellipse,draw,label=above:{$Y^*$},minimum height=3.5cm, minimum width=1cm] at (0,0) (Y) {};
  \node[circle,draw,fill,scale=.5,label=below:{$y_1$}] at (0,1.1) (y1) {};
  \node[circle,draw,fill,label=below:{$y_2$},scale=.5] at (0,.3) (y2) {};
  \node[circle,draw,fill,label=below:{$y_3$},scale=.5] at (0,-.3) (y3) {};
  \node[circle,draw,fill,label=below:{$y_4$},scale=.5] at (0,-1.1) (y4) {};
  \node[circle,draw,fill,label=above:{$o_3$},scale=.5] at (-1,1) (o4) {};
  \node[circle,draw,fill=blue,scale=.5,label=left:{$v_2$}] at (-3,.5) (x2) {};

  \draw (x2) [dashed,blue,line width=1.5] to (o4);
  \draw (x2) [dashed,blue,line width=1.5] to (y2);
  \draw (x2) [dashed,blue,line width=1.5] to (y3);
  \draw (x2) [dashed,blue,line width=1.5] to (y4);

  \draw (y4) [dashed,blue,line width=1.5, bend right=50] to (y1);
  \draw (o4) [dashed,blue,line width=1.5] to (y4);

  \begin{pgfonlayer}{background}
    \draw [line width=5, color=red!30](x2.center) to (y2.center);
    \draw [line width=5, color=green!30](x2.center) to (o4.center) to (y4.center);
    \draw [line width=5, color=olive!30](x2.center) to (y4.center) [bend right=50] to  (y1.center);
    \draw [line width=5, color=orange!30](x2.center) to (y3.center);
  \end{pgfonlayer}
  \node[circle,scale=1.3] at (-1.5,3) (m1) {$K(v_2)$};

  \node[circle,scale=4] at (2,0) (m1) {$\Rightarrow$};

  \node[ellipse,draw,label=above:{$Y^*$},minimum height=3.5cm, minimum width=1cm] at (7,0) (Y) {};
  \node[ellipse,draw,label=above:{$Y_C$},minimum height=3.5cm, minimum width=.5cm] at (9,0) (Y) {};
  \node[circle,draw,fill=red,scale=.5,label=right:{$f(v_1)$}] at (9,1) (fv1) {};
  \node[circle,draw,fill=blue,scale=.5,label=right:{$f(v_2)$}] at (9,0) (fv2) {};
  \node[circle,draw,fill=green,scale=.5,label=right:{$f(v_3)$}] at (9,-1) (fv3) {};
  \node[circle,draw,fill,scale=.5,label=below:{$y_1$}] at (7,1.1) (y1) {};
  \node[circle,draw,fill,label=below:{$y_2$},scale=.5] at (7,.3) (y2) {};
  \node[circle,draw,fill,label=below:{$y_3$},scale=.5] at (7,-.3) (y3) {};
  \node[circle,draw,fill,label=below:{$y_4$},scale=.5] at (7,-1.1) (y4) {};
  \node[circle,draw,fill,label=above:{$o_3$},scale=.5] at (6,1) (o4) {};
  \node[circle,draw,fill=blue,scale=.5,label=left:{$v_2$}] at (4,.5) (x2) {};

  \draw (x2) [dashed,blue,line width=1.5] to (o4);
  \draw (x2) [dashed,blue,line width=1.5] to (y2);
  \draw (x2) [dashed,blue,line width=1.5] to (y3);
  \draw (x2) [dashed,blue,line width=1.5] to (y4);

  \draw (y4) [dashed,blue,line width=1.5] to (fv2) to (y1);
  \draw (o4) [dashed,blue,line width=1.5] to (y4);

  \begin{pgfonlayer}{background}
    \draw [line width=5, color=red!30](x2.center) to (y2.center);
    \draw [line width=5, color=green!30](x2.center) to (o4.center) to (y4.center);
    \draw [line width=5, color=olive!30](x2.center) to (y4.center) to (fv2.center) to  (y1.center);
    \draw [line width=5, color=orange!30](x2.center) to (y3.center);
  \end{pgfonlayer}
  \node[circle,scale=1.3] at (7,3) (m1) {$K(v_2)$};
  }\\[0.6cm]

  \caption{Removing edges in $Y^*$ using $Y_C$.}
  \label{fig:removing edges in Y* using YC}
\end{figure}

Now, Claim~\ref{claim:if y' is small, then y^c is big} implies 
$|Y_C| > \ell - |\overline{X}_a| = |X^*|$.
Thus, there is an injection $f \colon X^* \to Y_C$, and we can replace in~$H$
every path $v u w \in K(v)$ for which $uw \in Y^*$ with the path $v u f(v) w$
(see Figure~\ref{fig:removing edges in Y* using YC}).
Notice that no edge is used more than once because~$E(Q(v))$ is a matching
in~$Y^*$.
Applying this replacement for every $v \in X^*$ yields that~$H$ contains an
immersion of $K_{X^*,Y^*}$ and~$G$ contains the desired immersion of 
$K_{\ell,\lceil n/2\rceil - \ell}$.
This concludes the proof of Theorem~\ref{theo:main}.

\section{Evidence for Conjecture \ref{conj:chromatic}}
\label{sec:consequences}

In this section we show that Theorem~\ref{theo:main} implies
Corollary~\ref{cor:maincorollary} (see~Corollary~\ref{cor:equiv}), and prove
Proposition~\ref{prop:ellequals2}.
For the first proof, we use arguments of Plummer \textit{et al.}~\cite{PlummerST2003}.
We need the following theorem and, for that, a definition: a graph~$G$ is said
to be \defi{$k$-critical} if $\chi(G) = k$ and $\chi(G-v) < k$, for every $v \in V(G)$.

\begin{theorem}[Stehl\'ik \cite{Stehlik2003}]\label{thm:critical}
  Let $G$ be a $k$-critical graph such that its complement is connected.
  For every $v \in V(G)$, $G-v$ has a $(k-1)$-coloring in which every color class
  contains at least two vertices.
\end{theorem}

\begin{corollary}
  \label{cor:equiv}
  Let $k,\ell$ be positive integers with $k \ge \ell+1$, and let $G$ be a
  $k$-chromatic graph with independence number~2.
  Then $G$ contains an immersion of $K_{\ell,k-\ell}$.
\end{corollary}
\begin{proof}
  Assume that $G$ is a counterexample that minimizes the number of
  vertices.
  It suffices to show that it is also a counterexample to Theorem~\ref{theo:main}.
  In order to apply Theorem~\ref{thm:critical}, we show that $G$ is $k$-critical
  and that its complement is connected.
  If $k = \ell+1$, we use that any graph with chromatic number~$k$ has 
  $K_{1,k-1} = K_{k-\ell,\ell}$ as a subgraph.
  Thus, we can assume that $k \geq \ell+2$.

  Let us first show that $G$ is $k$-critical.
  Indeed, if $G$ is not $k$-critical, then there would be some $v \in V(G)$
  such that $\chi(G-v) = k$.
  Then $G-v$ has at least $k$ vertices and is either a complete graph or has independence
  number~2, which, by choice of $G$, means it contains an immersion of
  $K_{\ell,k-\ell}$, a contradiction.

  Next, we show that the complement of $G$ is connected.
  If it were not, $V(G)$ could be partitioned into two sets $V_1, V_2$ such that
  $uv \in E(G)$ whenever $u \in V_1, v\in V_2$.
  Let $k_1, k_2$ be the chromatic numbers of $G[V_1]$ and $G[V_2]$, respectively,
  and note that $k = k_1+k_2$.
  As $k\geq \ell+2$, let $\ell_1,\ell_2$ be such that $\ell_1+\ell_2=\ell$,
  $k_1\geq \ell_1+1$, and $k_2\geq \ell_2+1$.
  Then, $G[V_1]$ and $G[V_2]$ must contain an immersion of
  $K_{\ell_1,k_1-\ell_1}$ and an immersion of $K_{\ell_2,k_2-\ell_2}$,
  respectively, because~$G$ is a counterexample that minimizes the number of
  vertices.
  But every vertex in $V_1$ is adjacent in $G$ to every vertex in $V_2$.
  Thus, $G$ must contain an immersion of
  $K_{\ell_1,k_1-\ell_1,\ell_2,k_2-\ell_2}$, which in turn implies an immersion
  of $K_{\ell_1+\ell_2,k_1-\ell_1+k_2-\ell_2}=K_{\ell,k-\ell}$.
  This contradicts the fact that $G$ is a counterexample.
  Therefore, the complement of $G$ is connected.

  Finally, by Theorem~\ref{thm:critical}, $G-v$ has a $(k - 1)$-coloring such
  that each color class has at least two vertices.
  But since~$G$ has independence number~2, each color class has size exactly~2.
  Thus we have $|V(G)| = 2k-1$.
  In particular, 
  $\bigl\lceil \frac{|V(G)|}2 \bigr\rceil = \frac {|V(G)|+1}2 = k$, meaning 
  that~$G$ does not contain an immersion of 
  $K_{\ell, \lceil |V(G)|/2 \rceil - \ell}$ and is a counterexample to Theorem~\ref{theo:main}, 
  as desired.
\end{proof}

Finally, we give a proof of Proposition~\ref{prop:ellequals2}, which claims
that every graph with chromatic number at least~3 contains an immersion of
$K_{1,1,\chi(G)-2}$.

\begin{proof}[Proof of Proposition~\ref{prop:ellequals2}.]
  Assume $\chi(G)=s+2 \geq 3$, and let $f \colon V(G)\to [s+2]$ be a proper
  coloring of~$G$ that minimizes the number of vertices colored with $s+2$.
  Let~$v$ be a vertex with $f(v)=s+2$.
  By the minimality condition, $v$ has at least one neighbor of each color.
  Furthermore, assume that~$f$ minimizes the number of vertices in~$N(v)$ colored
  with $s+1$.

  Let $vw \in E(G)$ with $f(w)=s+1$.
  For $1 \leq i \leq s$, consider the subgraph induced by vertices colored with~$i$
  or~$s+1$, and let $K_{wi}$ be the connected component containing vertex~$w$.
  If $K_{wi}$ does not contain neighbors of~$v$ colored with~$i$, then we can
  switch the colors of $K_{wi}$ (from $s+1$ to~$i$ and vice versa), obtaining a
  new coloring of~$G$ with less neighbors of~$v$ colored with $s+1$.
  Hence $K_{wi}$ must contain at least one neighbor of~$v$ colored with~$i$.
  Label one such vertex~$w_i$.

  For $1 \leq i \leq s$, let~$P_i$ be a path in $K_{wi}$ joining~$w$ to~$w_i$.
  Notice that if $1 \leq i \neq j \leq s$, then the path~$P_i$ is edge-disjoint
  from the path~$P_j$, since~$P_i$ uses edges with endpoints colored~$i$ and~$s+1$
  while~$P_j$ uses edges with endpoints colored~$j$ and~$s+1$.
  Further notice that these paths do not use any edge with~$v$ as one of its
  endpoints.
  Hence
  $$ F =  \bigcup_{i = 1}^s E(P_i)  \cup  \{vw_i \colon 1 \leq i \leq s + 1 \} $$
  induces an immersion of $K_{1, 1, s}$.
\end{proof}

\section{Immersions of matchings}
\label{sec:matching-immersions}

In this section, we discuss a related result that we obtained when proving
Theorem~\ref{theo:main}, and that we believe could have further applications.
Let~$G$ be a graph with a matching~$M$, and let $A,B \subseteq V(G)$ be
disjoint sets of vertices of~$G$.
We say that a matching~$M$ is an \emph{$(A,B)$--matching} if \(|M| = |A| = |B|\) and every edge in~$M$ has a
vertex in~$A$ and another vertex in~$B$.
We prove the following.

\begin{theorem}
\label{theo:matching-immersion}
  Let~$k$ and \(s\) be positive integers with \(s\leq k/2 + 1\), let $G = K_{2k}$, and let $(A_1,B_1), \ldots, (A_s,B_s)$
  be pairs of sets of vertices of~$G$ such that for every $i \in \{1,\ldots,s\}$ we have
  \begin{enumerate}[(a)] 
    \item $k \geq |A_i| = |B_i| \geq 2(i-1)$; and
    \item $A_i \cap B_i = \emptyset$.
  \end{enumerate}
  Then~$G$ contains edge-disjoint matchings $M_1,\ldots, M_s$ such that~$M_i$
  is an $(A_i,B_i)$--matching for $i \in \{1,\ldots,s\}$.
\end{theorem}
\begin{proof}
  Assume $M_1,\ldots,M_{i-1}$ were found.
  Let $U = M_1 \cup \cdots \cup M_{i-1}$ be the set of used edges and let
  $F = E(G[A_i \cup B_i]) \setminus U$ be the set of free edges.
  We use Hall's Theorem.
  For that, let $X\subseteq A_i$.
  We prove that $|N_F(X) \cap B_i| \geq |X|$.
  Suppose, for a contradiction, that $|N_F(X) \cap B_i| < |X|$.
  Since $d_U(u) \leq i-1$ for every $u \in V(G)$, we have 
  $|N_F(u) \cap B_i| \geq |B_i|-(i-1)$. By choosing \(u\in X\) we obtain
  \begin{equation}
  \label{eq:matching-immersions:eq1}
    |X| > |N_F(X) \cap B_i| \geq |N_F(u)\cap B_i| \geq |B_i| - (i-1) \, .
  \end{equation} 
  On the other hand, if $u \in B_i \setminus N_F(X)$, then 
  \begin{equation*}
  \label{eq:matching-immersions:eq2}
    |A_i| - |X| \geq d_F(u) \geq |A_i| - (i-1) \, ,
  \end{equation*} 
  from where we get $|X| \leq i-1$,
  which, together with~\eqref{eq:matching-immersions:eq1}, gives $|B_i| < 2(i-1)$, 
  a contradiction.
\end{proof}

As we conclude this paper, we raise the following question, which if true, would 
be a generalization of Theorem~\ref{theo:matching-immersion}.

\begin{question*}
Let~$k$ and \(s\) be positive integers with \(s\leq k\), let $G = K_{2k}$, and let $(A_1,B_1), \ldots, (A_s,B_s)$
  be pairs of sets of vertices of~$G$ such that for every $i \in \{1,\ldots,s\}$ we have $|A_i| = |B_i| \leq k$ and $A_i \cap B_i = \emptyset$.
  Does~$G$ contain edge-disjoint immersions $M_1,\ldots, M_s$ such that $M_i$ is an immersion
  of an $(A_i,B_i)$--matching for $i \in \{1,\ldots,s\}$?
\end{question*}

\section*{Acknowledgement}

We thank the referees for their valuable comments and suggestions.
We thank Feri Kardo\v{s}, Matías Pavez-Signé, Vinicius dos Santos, and Maya
Stein for many stimulating discussions.

\bibliographystyle{amsplain}
\bibliography{bibfile}

\end{document}